\documentclass[a4paper,11pt]{article} 

\usepackage{fullpage}         
\usepackage[english]{babel}
\usepackage[utf8]{inputenc}
\usepackage[style=alphabetic, backend=bibtex8, natbib]{biblatex}
\usepackage{authblk}          
\usepackage{amssymb}          
\usepackage{amsthm}	          
\usepackage{amsmath}
\usepackage{upgreek}          
\usepackage{mathrsfs}
\usepackage{graphicx}         
\usepackage{float}            
\usepackage{enumitem}
\usepackage{xcolor}           
\usepackage{hyperref}
\usepackage{algpseudocode}
\usepackage{algorithm}
\usepackage[title]{appendix}

\usepackage{caption}
\usepackage{subcaption}
\hypersetup{
    linktocpage = true,
    colorlinks  = true,                        
    linkcolor   = blue,
    citecolor   = blue,
    urlcolor    = blue
}

\numberwithin{equation}{section}

\DeclareMathOperator*{\vspan}{span}
\DeclareMathOperator*{\tr}{tr}
\DeclareMathOperator*{\ran}{ran}

\DeclareMathOperator*{\kl}{KL}

\theoremstyle{definition}
\newtheorem{definition}{Definition}[section]

\theoremstyle{remark}

\newtheorem{remark}[definition]{Remark}

\theoremstyle{plain}
\newtheorem{lemma}[definition]{Lemma}
\newtheorem{corollary}[definition]{Corollary}
\newtheorem{theorem}[definition]{Theorem}
\newtheorem{proposition}[definition]{Proposition}

\begin{document}
\title{Posterior contraction rates for banded approximations of Gaussian processes\footnote{This paper has been accepted for publication at the NeurIPS 2026 conference.}}
\renewcommand\Affilfont{\itshape\small}
\author[1]{Omar Kassi}
\affil[1]{Department of Decision Sciences, Bocconi University,
          Bocconi Institute for Data Science and Analytics (BIDSA), 
          \normalfont \texttt{omar.kassi@unibocconi.it}
         }
\author[2]{Bernhard Stankewitz}
\affil[2]{Department of Mathematics, University of Potsdam,
          \normalfont \texttt{bstankewitz@posteo.de}
         }
\author[3]{Botond Szab{\'o}}
\affil[3]{Department of Decision Sciences, Bocconi University,
          Bocconi Institute for Data Science and Analytics (BIDSA), 
          \normalfont \texttt{botond.szabo@unibocconi.it}
         }
\date{\today}
\maketitle
\abstract{
  Gaussian process regression models are widely used in modern statistics and machine learning due to their flexibility, interpretability, built-in uncertainty quantification, and strong theoretical foundations. However, training Gaussian processes (GPs) has a computational cost of $O(n^3)$, while the memory and prediction costs are also $O(n^2)$. 
  For large datasets commonly encountered in practice, this becomes infeasible and necessitates the use of approximate posterior methods.
  In this work, we propose a class of approximations based on sparse covariance matrices, where the kernel matrix is sparsified according to the spatial proximity of the design points. 
  We analyze the accuracy of the resulting approximate posterior distribution and show that the proposed approach preserves optimal statistical inference properties, provided that the dependencies among a sufficient number of neighboring points for each feature are retained. 
  We derive explicit sufficient and necessary conditions on the number of neighbors required in a general setting and apply these results to standard kernels, including the squared exponential and Mat\'ern kernels. Our theoretical findings are supported by numerical experiments.
}


\section{Introduction}
\label{sec_Introduction}

Gaussian Processes (GPs) arguably form one of the standard and most popular Bayesian non-parametric approaches for regression and classification tasks. 
They have a broad spectrum of applications, ranging from disease risk mapping in the Malaria Atlas Project \cite{BhattEtal2017_ImprovedPredictionAccuracy} through understanding properties of molecules and chemicals \cite{DeringerEtal2021_GPRegression} to aircraft designs \cite{ArnsonEtal2025_PredictingAircraftDesign}. 
Their high popularity stems from the built-in uncertainty quantification, their effectiveness even with small data sizes, their flexibility, and improved interpretability compared to black box procedures.

In this paper, we focus on the random design nonparametric regression model, in which we observe pairs of random variables $(X_i,Y_i)_{i=1,...,n}$ satisfying
\begin{equation}
  \label{eq_Intro_GPRegressionModel}
  Y_i = f(X_i) + \varepsilon_i, \quad \text{with} \quad 
                 \varepsilon_i \sim \mathcal{N}(0, \sigma^2) \text{ i.i.d., }
                           X_i \sim G  \text{ i.i.d., }
                 \qquad  i = 1, \dots, n,
\end{equation}
where $f : \mathcal{X} \to \mathbb{R}$, $\mathcal{X}\subset\mathbb{R}^d$, denotes the unknown regression function of interest and $G$ denotes the distribution of the each design point.
Also \( \varepsilon = (\varepsilon_i)_{i \le n} \) and \( X = (X_i)_{i \le n} \) are independent.
In the Bayesian approach, $f$ is endowed with a flexible nonparametric prior.
Gaussian Processes provide a convenient and highly efficient choice for the prior measure.
Due to conjugacy, the corresponding posterior is also a Gaussian Process and can be explicitly derived, see e.g. \cite{RasmussenEtal2006_GaussianProcesses}.
Denoting by $k(\cdot, \cdot)$ the covariance kernel of a centered GP used as the prior, the mean and covariance functions of the corresponding posterior take the form
\begin{equation} \label{eq_Posterior}
\begin{split}
  x       & \mapsto k(X,x)^{\top} (K + \sigma^2 I_n)^{-1} Y, \qquad \qquad \qquad \qquad \qquad \ \ x \in \mathcal X,  \\
  (x, x') & \mapsto k(x,x') - k(X,x)^{\top} (K + \sigma^2 I_n)^{-1} k(X,x'), \qquad x, x' \in \mathcal X,
\end{split}
\end{equation}
where $ K := ( k(X_i, X_j) )_{i, j = 1,..,n} \in \mathbb{R}^{n \times n}$ denotes the empirical prior covariance matrix, $ I_n $ the $ n $-dimensional identity matrix, $ Y = (Y_1, \dots,Y_n)^{\top} \in \mathbb{R}^n $ the observation vector and $ k(X, x) = ( k(X_1, x), \dots, k(X_n,x) )^{\top} \in \mathbb{R}^n$.
Note that the inversion of the matrix $ K_{\sigma}: = K + \sigma^2 I_n \in \mathbb{R}^{n \times n}$ has computational complexity $O(n^3)$.
Hence the computation of the GP posterior has training cost of $O(n^3)$, memory cost $O(n^2)$ and prediction cost $O(n^2)$.
These render the computation infeasible or at least undesirable for large sample sizes.
Since in modern applications the number of available data can be extensive this hinders the practical applicability of GPs.

To overcome this issue often scalable alternatives are considered, where the GP posterior (or prior) is approximated with a low rank/sparse process. 
There is a wide range and ever increasing number of such approaches, including distributed GPs \cite{JacobsEtal1991_AdaptiveMixturesLocalExperts, GramacyLee2008_BayesianTreedGaussianProcess, Tresp2000_GeneralizedBayesianCommitteeMachine, ScottEtal2016_BayesBigDataConsensusMonteCarlo}, variational Bayes methods using inducing variables \cite{Titsias2009_InducingVariablesVB}, covariance tapering \citep{FurrerEtal2006_CovarianceTaperingInterpolation, KaufmanEtal2008_CovarianceTaperingLikelihoodEstimation, Stein2013_StatisticalPropertiesCovarianceTapers}, and Vecchia approximation \cite{Vecchia1988_EstimationModelIdentificationSpatialProcesses, DattaEtal2016_HierarchicalNearestNeighborGaussianProcesses, KatzfussEtal2020_VecchiaApproximationsGaussianProcessPredictions}. 
These methods typically rely on carefully tuning various sparsity, low-rank, and localization parameters, trading off accuracy against computational complexity. 
The theoretical understanding of these approaches and the optimal tuning has been investigated in recent years, e.g. for distributed GPs \cite{SzaboVanZanten2019_AsymptoticAnalysisDistributedNonparametricMethods, GuhaniyogiEtal2022_DistributedBayesianVaryingCoefficientModeling, SzaboEtal2025_AdaptationSpatiallyDistributedGaussianProcesses}, Vecchia GPs \cite{SzaboZhu2026_VecchiaGaussianProcesses} variational GPs \cite{BurtEtal2019_RatesConvergenceSparseVariationalGP, NiemanEtal2022_ContractionRates}, but our fundamental understanding is still limited.

In our paper, we consider a novel approximation method inducing sparsity to the prior covariance matrix $ K $ and hence substantially speeding up the inversion of the matrix and as a consequence reducing the cost of training, prediction and storage.\footnote{While writing up our manuscript we learned about the recently published conference paper \cite{EhrhardtTobar2026_AcceleratedTrainingGaussianProcesses} proposing a similar method in case of squared exponential covariance kernel. The authors investigate some practical aspects of this procedure, but no statistical guarantees were derived and extension to other GP kernels were also not considered.}
To preserve the statistical properties of the posterior and minimize information loss, for each feature $X_i$, we retain its dependence on its closest neighbors within a radius \( r > 0 \) while setting all other entries of the covariance matrix \( K \) to zero.
The idea of preserving the information available in the closest neighbours while discarding the less relevant peripheral data has been widely used in other approximation methods in various manners. 
For instance, in the nearest neighbour Vecchia approximation, the precision matrix of the prior is sparsified in a similar way.
In covariance tapering, the covariance matrix is multiplied by a sparse, compactly supported correlation function, forcing small, distant covariances to zero while maintaining a positive definite matrix. 
While our approach is closely connected to the latter one, it has substantial algorithmic and conceptual differences, see the next section discussing them in detail.

Our paper investigates both the methodological and theoretical aspects of the proposed approximation approach. 
We derive general theoretical guarantees and limitations of the method for statistical inference, focusing on recovering the functional parameter of interest in the regression model and on the approximation accuracy of the posterior. 
We connect our results to the general theory of nonparametric Bayesian methods using Gaussian Process priors and obtain sufficient and necessary conditions for the approximation to maintain the same level of accuracy as the original posterior while substantially reducing the computational complexity. 
Our general, abstract results are then applied to standard covariance kernels revealing the trade-off between computational speed and approximation accuracy in concrete examples. 
Our theoretical results then can be used as a guideline in practice for tuning the approximations in an optimal manner instead of following an ad hoc empirical approach.     

\textbf{Organization}.
In Section \ref{sec_Algorithms} we introduce the approximation approach rigorously and discuss its methodological aspects. 
We provide general theoretical guarantees for approximation accuracy and statistical inference in Section \ref{sec_MainResults} and apply these results for the popular squared exponential and Mat\'ern kernels in Section \ref{sec_Examples}. 
Our theoretical findings are demonstrated with numerical examples considering synthetic data sets. 
The proof of the results together with additional examples and numerical studies are deferred to the Appendices.

\textbf{Notations}.
For two positive sequences $(a_n)$ and $(b_n)$, we write $a_n \lesssim b_n$ if there exists a constant $C>0$, independent of $n$, such that $a_n \leq C b_n$ for all $n$. 
We further write $a_n \asymp b_n$ if both $a_n \lesssim b_n$ and $b_n \lesssim a_n$ hold. 
When needed, we explicitly use constants $c, C > 0$ that do not depend on $n$ and may vary from line to line. 
For a function $f \in L^2(G)$, its evaluation at the design points $(X_1,\dots,X_n)$ is denoted by $\mathbf{f} := (f(X_1),\dots,f(X_n)) \in \mathbb{R}^n$. 
Furthermore, we denote by $K=k(X,X)\in\mathbb{R}^{n\times n}$ the empirical covariance matrix and by $K_\sigma=K+\sigma^2 I_n\in\mathbb{R}^{n\times n}$ the matrix appearing the GP posterior. 
For two probability distributions $\mathbb P$, $\mathbb Q$,
$ \kl(\mathbb P, \mathbb Q) = \int \log ({\rm d}\mathbb P/ {\rm d}\mathbb Q) d \mathbb P$ denotes the Kullback-Leibler (KL) divergence.
Given a probability measure $\mathbb P$, the distribution of a random variable $X$ under $\mathbb P$ is denoted by $\mathbb P^X$ and $\mathbb E_X$ denote the associated expectation. 
Likewise, for another random variable $Y$, the conditional distributions of $Y\mid X$ and $Y\mid X = x$ are denoted by $\mathbb P^{Y \mid X}$ and $\mathbb P^{Y \mid X = x}$, respectively. 


\section{Algorithms}
\label{sec_Algorithms}

For many commonly used kernels, correlations between observations decay rapidly as the distance between design points increases. 
This suggests that entries of the kernel matrix corresponding to distant points contribute little to the posterior and may be neglected. 
Therefore, a straightforward idea to reduce the computational complexity of GP regression is to sparsify the kernel matrix \( K \) by setting small correlations to zero when design points are far apart. 
This leads to the following approximation \( \widehat K = (\widehat K)_{i, j \le n} \) of \( K \) with 
\begin{align}
  \label{eq_DefKHat}
  \widehat{K}_{ij} = \begin{cases}
                       k(X_i, X_j), & \| X_i - X_j \| \le r, \\
                       0, & \text{otherwise},
                     \end{cases}
\end{align}
where \( r > 0 \) is a fixed radius.
In dimension \( d = 1 \), when the design points can be assumed to be ordered on the real line, this construction reduces to a banded matrix with maximal bandwidth \( m = m_r \) depending on \( r \). 
In higher dimensions, where no natural ordering is available, the approach more closely related to covariance tapering methods \cite{FurrerEtal2006_CovarianceTaperingInterpolation}.
Covariance tapering typically constructs a new kernel by multiplying the original kernel with a compactly supported function, and then performs full Bayesian inference under this modified prior. 
In contrast, our perspective retains the original kernel as the underlying model and views \( \widehat{K} \)  as a computational approximation of \( K \). 
As a consequence, the resulting posterior should be interpreted as an approximation to the true posterior rather than the exact posterior under a modified prior. 
This viewpoint is particularly useful when analyzing how the sparsity parameter, i.e. the radius \( r \), should scale with the sample size.
Defining $\widehat{K}_\sigma := \widehat{K} + \sigma^2 I_n$, we obtain the approximate posterior distribution $\Psi_{r}$, whose mean and covariance functions are given by
\begin{equation}
  \label{eq_ApproximatePosterior}
  x \mapsto k(X,x)^{\top}\widehat{K}_\sigma^{-1}Y, \qquad 
  (x,x') \mapsto k(x,x') - k(X,x)^{\top}\widehat{K}_\sigma^{-1}k(X,x'), 
  \qquad x, x' \in \mathcal X.
\end{equation}
\( \Psi_r \) is a well defined posterior as long as \( \| \widehat K - K \|_{\text{op}} < \sigma^2 \), see Appendix \ref{ssec_WellDefinednessOfTheApproximatePosterior}.
Setting, analogously to the one-dimensional case,  \( m = m_r \) equal to the maximal number of non-zero entries in a row of \( \widehat K \), storing \( \widehat{K}_\sigma \) scales as \( \mathcal{O}(mn) \) in terms of memory compared to \( \mathcal{O}(n^2) \) for the full kernel matrix.
Importantly, the full matrix \( K \) does not need to be stored explicitly; only local kernel evaluations are required to construct \( \widehat{K} \) and to evaluate terms such as \( k(X,x) \).

The main advantage of the sparsified matrix \( \widehat{K}_\sigma \) is that it contains at most \( \mathcal{O}(mn) \) nonzero entries. 
This structure is particularly well-suited for iterative linear solvers such as the conjugate gradient (CG) algorithm. 
Each iteration of CG only requires matrix-vector multiplications of the form \( \widehat{K}_\sigma u \), along with vector operations and inner products. 
When \( \widehat{K}_\sigma \) is dense, such multiplications cost \( \mathcal{O}(n^2) \) operations. 
In contrast, the sparsity of \( \widehat{K}_\sigma \) reduces this cost to \( \mathcal{O}(mn) \), since only the nonzero entries contribute. 
Then, the total computational cost of solving the system is therefore \( \mathcal{O}(mn \cdot N_{\mathrm{iter}}), \) where \( N_{\mathrm{iter}} \) is the number of iterations required for convergence.  

A key question is how to choose the parameter \( r > 0 \).
Larger values of \( r \) improve the approximation of the true kernel matrix but also increase computational cost, while smaller values enhance scalability at the expense of approximation accuracy.  
In particular, one expects that \( m = m_r \) should grow with \( n \) to ensure that the approximation error vanishes asymptotically, while remaining as small as possible to retain computational gains. 
This interplay between \( r \), \( m \) and \( n \) will be reflected in the theoretical results developed in subsequent sections and analyzed also numerically in a range of examples.


\section{Main results}
\label{sec_MainResults}

We now provide theoretical guarantees for the approximation method described in the previous section. 
In our theoretical analysis, we consider the asymptotic regime where the sample size $n$ tends to infinity. 
The finite sample behaviour of the procedure is investigated numerically in the next sections. 
We derive how well the proposed scalable alternative approximates the original posterior distribution and how much it preserves its statistical inference properties. 
To analyze the accuracy of the method, we consider the standard frequentist framework where we assume that there exists a deterministic true regression function $f_0$ in the model from Equation \eqref{eq_Intro_GPRegressionModel}. 
Then, we derive how fast the approximate posterior recovers the underlying function $f_0$ as the sample size increases. 
This is measured via the speed of contraction of the approximate posterior around $f_0$. 

First, we consider the approximation accuracy of the proposed scalable alternative $\Psi_{r}$ to the GP posterior $\Pi_n(\cdot | X, Y)$.  
The geometry of the prior Gaussian Process $W$ is described by the Reproducing Kernel Hilbert Space (RKHS) $\mathbb H$ associated to its covariance kernel $k$. 
Indeed, the RKHS \( \mathbb{H} \) of \( W \) is obtained simply as the completion of the linear span of functions \( \{ k(x, \cdot): x \in \mathcal{X} \} \) endowed with the inner product 
\begin{align}
    \Big\langle \sum_{i = 1}^p a_i k(t_i, \cdot), \sum_{j = 1}^q b_j k(s_j, \cdot)
    \Big\rangle_{\mathbb H}
  =
    \sum_{i = 1}^p \sum_{j = 1}^q a_i b_j k(t_i, s_j),
\end{align}
for arbitrary  $a\in \mathbb{R}^{p}, b \in \mathbb{R}^{q}, t_1,..,t_p, s_1,...,s_q\in\mathcal{X}$.
The associated RKHS norm is then defined as $\|g\|_{\mathbb H}^2=\langle g,g\rangle_{\mathbb{H}}$ for $g\in\mathbb{H}$. 
For a detailed introduction to RKHSs and their analytical properties we refer the reader to Section \ref{ssec_ReproducingKernelHilbertSpaces} in the Appendix. 
Then, the statistical properties of the posterior corresponding to the GP prior is encoded in the behaviour of the so called concentration function
\begin{equation}
  \label{eq_Main_ConcentrationFunction}
  \varphi_{f_0}(\eta) = \inf_{h\in\mathbb H:\,\|h-f_0\|<\eta}\|h\|_{\mathbb H}^2
                      - \log \mathbb{P}\big(\|W\| < \eta\big),
\end{equation}
see  \cite{vdVaartvZanten2008Rate}. The first term determines the approximation accuracy of the function $f_0$ with the RKHS. This term can be associated to the bias or approximation error of the procedure. The second term in the previous display is the log small ball probability.
This term can be associated to the spread of the Bayesian procedure.
Then, in view of \cite{vdVaartvZanten2008Rate}, the posterior corresponding to the GP prior contracts around the true parameter $ f_0 $ with rate $ \eta_n $, if the inequality 
\begin{equation}
 \label{eq_Main_ConcentrationFunctionInequality}
  \varphi_{f_0}(\eta_n) \le C_\varphi n \eta_n^2
\end{equation}
holds for some constant \( C_\varphi > 0 \).
Typically the optimal contraction rate can be expressed as the fastest rate \( \eta_{n} \) that satisfies the previous inequality. 
This is a standard condition describing the statistical inference properties of the GP posterior and hence is a natural assumption in our case as well.

Next, we express the approximation error of the proposed scalable approach of the original GP posterior with respect to the estimation error between the sparsified covariance matrix $\widehat{K}$ and the true covariance matrix $K$. 
The latter quantity naturally appears in the approximation error and can be controlled based on the spectral and spatial properties of the prior covariance kernel $k$.

\begin{theorem}[Kullback-Leibler bound]
	\label{thm_KullbackLeiblerBound}
	Assume that the concentration function satisfies condition \eqref{eq_Main_ConcentrationFunctionInequality} for some sequence $\eta_n$. 
  Moreover, suppose that the threshold radius \( r = r_n \) is chosen such that
	\begin{equation}
    \label{eq_Main_PerturbationCondition}
		\| \widehat{K} - K \|_{\mathrm{op}} \leq \overline c < \sigma^2/2
	\end{equation}
	with $\mathbb P^{\otimes n}_{f_0}$-probability converging to one. 
  Then there exists a constant $\mathfrak C>0$ depending on $\sigma^2, d$, and $\overline c$ such that, for any sequence $M_n \to \infty$, the approximate posterior $\Psi_{r_n}$ satisfies
  \begin{align}
        \kl\!\left(\Psi_{r_n}, \Pi_n(\cdot \mid X, Y)\right)
    \le
        \mathfrak C\, M_n \| \widehat{K} - K \|_{\mathrm{op}}^2\,
                          \left\{ \tr(KK_\sigma^{-1})+ n \eta_n^2 \right\},
  \end{align}
	with probability converging to one under $\mathbb P^{\otimes n}_{f_0}$ as $n \to \infty$.  
\end{theorem}

\noindent The above theorem, which is proven in Appendix \ref{prf_KullbackLeiblerBound}, shows that in case that $ \| \widehat{K} - K \|_{\mathrm{op}}^2 = o\big( (\tr(KK_\sigma^{-1})+ n \eta_n^2)^{-1}\big)$ with probability converging to one, the Kullback-Leibler divergence  between the approximate and true posteriors tends to zero.
The operator norm $\| \widehat{K} - K \|_{\mathrm{op}}$ decreases as the threshhold radius $ r > 0 $ in the approximation $\widehat{K}$ increases. 
The exact decay depends both on the spatial and spectral properties of the covariance kernel $k$. 
We will investigate this interplay in case of several examples, including both polynomially and exponentially decaying eigenvalues and various spatial behavior of the covariance kernel. 

Controlling the Kullback-Leibler approximation error of the posterior provides us the means to derive similar contraction rate results for the approximation $\Psi_{r_n}$ as the original posterior $\Pi_n(\cdot| X, Y)$, see Theorem 5 of \cite{RaySzabo2022_VariationalBayesHighDimensionalRegression}. 
Formally, the sequence \( (\Psi_{r_n})_{n \in \mathbb{N}} \) of distributions on \( L^{2}(G) \) contracts around the truth \( f_{0} \) with respect to a metric \( d \) if, for an arbitrary sequence \( M_n \to \infty \),
\begin{align}
  \Psi_{r_n}( d(\cdot, f_{0}) \ge M_{n} \eta_{n}) \to 0
  \qquad 
  \text{in  probability under } \mathbb{P}_{f_{0}}^{\otimes n} \text{ as } n \to \infty.
\end{align}
Intuitively, this means that asymptotically, \( \Psi_{r_n} \) does not put any mass outside of a ball of size \( O(\eta_{n}) \). 
In the following, we will consider the metric \( d = d_{H} \) given by 
\begin{align}
  d_{H}(f, g)^2 = \int (\sqrt{p_{f}} - \sqrt{p_{g}})^{2} d G \otimes \lambda, 
  \qquad f, g \in L^{2}(G),
\end{align}
which is induced by the Hellinger distance between the two Gaussian densities \( p_{f} \) and \( p_{g} \) with \( p_{f}(x, y) = \exp(- (y - f(x))^{2} / (2 \sigma^{2})) / \sqrt{2 \pi \sigma^{2}} \).
The next theorem provides contraction rates for the approximate posterior $\Psi_{r_n}$.
Its proof is in Appendix \ref{prf_ContractionRates}.

\begin{theorem}[Contraction rates]
  \label{thm_ContractionRates}
  In the situation of Theorem \ref{thm_KullbackLeiblerBound}, assume that additionally \( \tr K K_ \sigma^{-1} = O_{\mathbb P_{f_0}}(n \eta_n^2) \).
  Then, for any sequence \( M_n \to \infty \), the approximate posterior satisfies
  \begin{align}
    \Psi_{r_n}(d_{\text{H}}(\cdot, f_0) \ge M_n \eta_n) \to 0
  \end{align}
  in probability under \( \mathbb P_{f_0}^{\otimes n} \) as \( n \to \infty \).
\end{theorem}

\noindent Note that the condition $\tr(K K_\sigma^{-1}) = O_{\mathbb{P}_{f_0}}( n \eta_n^2)$ ensures that the Kullback-Leibler divergence between the approximation and the true posterior is of the order $O_{\mathbb{P}_{f_0}}( \| \widehat{K} - K \|_{\mathrm{op}}^2 n \eta_n^2)$. 
Then, under the assumption $ \| \widehat{K} - K \|_{\mathrm{op}} < \sigma^2 / 2$, the Kullback-Leibler divergence is bounded by a multiple of $ n \eta_n^2$, which in view of  \cite{RaySzabo2022_VariationalBayesHighDimensionalRegression}, is sufficient to derive the contraction rate $\eta_n$ for $\Psi_{r_n}$. 
Note that to obtain this rate, it is not required that the approximation error tends to zero in Kullback-Leibler divergence.
In fact, it can diverge to infinity. 
We discuss this disparity phenomenon in our examples.
In the next section, we demonstrate the above general results on concrete examples, providing both theoretical and numerical justifications.


\section{Examples}
\label{sec_Examples}

In this section, we present two explicit examples and provide a guide for selecting the neighbourhood size $ r $ in the approximation. 
The derived bounds depend on the smoothness of the underlying regression function $f_0$, the regularity of the GP, the spatial regularity of the covariance kernel and the dimension $d$ of the process. 
To make this precise, we recall two classical smoothness classes. 
The Hölder space $C^\alpha(\mathbb{R}^d)$, for $\alpha > 0$, consists of functions whose partial derivatives up to order $\lfloor \alpha \rfloor$ exist and are bounded, and whose derivatives of order $\lfloor \alpha \rfloor$ are Hölder continuous with exponent $\alpha - \lfloor \alpha \rfloor$.
The Sobolev space $H^\alpha(\mathcal X)$ is defined as the class of restrictions $h|_{\mathcal X}$ of functions $h : \mathbb{R}^d \to \mathbb{R}$ whose Fourier transform $\hat{h}$ satisfies
\[
\int_{\mathbb{R}^d} (1 + \|y\|^2)^\alpha |\hat{h}(y)|^2 \, d y < \infty.
\]
It is well known that, in contrast to regular parametric models, in a nonparametric setting the choice of the prior plays a crucial role in the asymptotic behaviour of the posterior. 
Centered Gaussian processes described by their covariance kernel typically involve one or more tuning parameters. 
These parameters can determine the smoothness of the process (e.g., in the case of Mat\'ern, Riemann–Liouville, or fractional Brownian motion kernels) or rescale the temporal or spatial behavior of the process. 
For example, set
$$
  W^{\tau_n}(x) = W(x/\tau_n),\quad x \in \mathcal X,
$$
where $\tau_n$ is a time re-scaling parameter possibly depending on the sample size. 
Then, one can tune the scale parameter $\tau_n$ to balance the bias and variance terms of the process, thereby preventing over-  and underfitting.

Below, we consider the rescaled version of two standard GPs, the squared exponential and Mat\'ern processes. 
We investigate the approximation error of the corresponding posteriors based on our approach, derive contraction rate results and study numerically the small to moderate sample size behaviour of the approximation procedure.

\subsection{The squared exponential kernel}
\label{ssec_TheSquaredExponentialKernel}

In this section, we focus on the rescaled squared exponential (SE) kernel in dimension $d$:
\begin{align}
  k_\tau(x, x') = \exp \big( - \tau^{- 2} \| x' - x \|^2 \big), \qquad x, x' \in \mathbb R^d,
\end{align}
where \( \| \cdot \| \) denotes the Euclidean norm.
For $\tau = 1$ the sample paths from the SE kernel are almost surely analytical.
Infinitely smooth GPs are not suitable to approximate a less regular regression function $f_0\in C^{\alpha}$, $\alpha>0$, as the bias of the procedure will be sub-optimally large. 
Therefore, we consider a rescaled version of the SE GP with scale parameter taken to be 
\( \tau = \tau_n =n^{-1 / (2 \alpha + d)} \).
This optimal choice  will shrink the sample paths of the SE GP  and allows for estimating $\alpha$-smooth functions $f_0$, see for instance \cite{vdVaartvZanten2007RescaledGPPriors} for 1-dimensional case and Lemma 8 of \cite{NiemanEtal2022_ContractionRates} for the d-dimensional case.
We work with this optimally rescaled process and approximate it via the proposed nearest neighbour technique. We show that, by retaining a sufficiently large number of neighbours, our approximation is accurate for the posterior with the KL divergence tending to zero.
Furthermore, the procedure preserves optimal statistical performance: the resulting posterior achieves minimax-optimal contraction up to a slowly varying term.

\begin{proposition}[Result for the squared exponential kernel]
  \label{prp_ResultForTheSquardExponentialKernel}
  In model \eqref{eq_Intro_GPRegressionModel}, assume that the design distribution \( G \) has a centered Gaussian Lebesgue density \( g \) on \( \mathcal X= \mathbb R^d \) and that \( f_0 \in C^\alpha(\mathbb R^d) \cap L^2(\mathbb R^d) \) for some \( \alpha > d/2 \).
  Then for the squared exponential kernel \( k_\tau \) with \( \tau = \tau_n = n^{-1 / (2 \alpha + d)} \), the approximate posterior with truncation radius \( r = r_n = \sqrt{2 \tau_n^2 \log n} \) satisfies that for arbitrary \( M_n \to \infty \), both
  \begin{align}
    \Psi_{r_n}(d_{\text{H}}(\cdot, f_0) \ge M_n n^{-\alpha / (d + 2 \alpha)} \log^{2 d} n) \to 0
    \quad \text{ and } \quad 
    \kl(\Psi_{r_n}, \Pi_n(\cdot | X, Y) \to 0
  \end{align}
  in probability with respect to \( \mathbb P_{f_0}^{\otimes n} \) for \( n \to \infty \).
  Further, the number of non-zero entries per row of \( \widehat K \) is at most of order \( m_n = n r_n^d \asymp n^{2 \alpha / (2 \alpha + d)} \log^{d/2} n\) with probability converging to one.
\end{proposition}

\noindent Proposition \ref{prp_ResultForTheSquardExponentialKernel}, which is proven in Appendix \ref{prf_ResultForTheSquardExponentialKernel}, shows that if the truncation radius is as above then the approximate posterior attains the optimal contraction rate, and the Kullback-Leibler divergence between the approximate and true posteriors converges to zero, i.e., the approximate posterior shares the essential statistical properties of the true posterior.
Regarding the computational complexity, the result guarantees that any row of the matrix \( \widehat K \) contains at most of order \( m_n = n r_n^d  \) non-zero entries.
In the context of the discussion in Section \ref{sec_Algorithms}, this reduces storage to \( O(m_n n) \) and the computational cost to \( O(m_n n N_{\text{iter}}) \).

Recall from Section \ref{sec_MainResults} that to obtain the same contraction rate as the original posterior one does not require that the Kullback-Leibler divergence between the two measures vanishes. 
In the present example, however, the gap between the optimal choices of $ r $ vanishes (up to log factors), see also Corollary \ref{cor_LowerBoundForTheSquaredExponentialKernel}.
This is due to the squared exponential tail decay of the covariance kernel $r_{\tau_n}(\|t\|)$ resulting in a substantial increase in the accuracy of the covariance matrix approximation $\widehat{K}$ in operator norm when increasing $ r $. 
Furthermore, note that as $\alpha$ increases to infinity, the exponent $2\alpha/(2\alpha + d)$ approaches 1, prohibiting any meaningful banding. 
This reflects the fact that smoother functions exhibit stronger long-range dependence, and thus require incorporating more points to be accurately approximated. 
On the other hand, as the dimension $d$ increases, the exponent decreases, which indicates that fewer neighbors are needed. 
Next we show that the required size of the truncation radius in Proposition \ref{prp_ResultForTheSquardExponentialKernel} is sharp.

\begin{corollary}[Lower bound for the squared exponential kernel]
  \label{cor_LowerBoundForTheSquaredExponentialKernel}
  In model \eqref{eq_Intro_GPRegressionModel}, assume that the design distribution \( G \) is standard Gaussian with density \( g \) on \( \mathcal X = \mathbb R^d \) and let \( \alpha > d/2 + 1 \).
  Then for the squared exponential kernel \( k_\tau \) with \( \tau = \tau_n = n^{-1 / (2 \alpha + d)} \), the truncated kernel matrix \( \widehat K \) with truncation radius \( r \) smaller than \( r_n = \tau_n /  \log^{1/(2d)} n \) satisfies that 
  \begin{align}
    \inf_{r < r_n} \sup_{f_0 \in C^{\alpha}(\mathbb R^d)}
    \| \mathbb{E}_{\varepsilon} K \widehat{K}_{\sigma}^{-1} Y - f_{0} \|_{n}
    \gtrsim \sqrt{\log n}\, n^{-\frac{d/2}{2\alpha+d}}
  \end{align}
  with probability converging to one under \( \mathbb P_{f_0}^{\otimes n} \).
\end{corollary}

\noindent The above corollary states that choosing even a logarithmic multiplicative factor less neighbours in the approximation, the posterior mean has sub-optimally large bias $ n^{-d / (2 (2 \alpha + d))} \gg n^{- \alpha / (2 \alpha + d)}$, for $\alpha > d/2 + 1$, with respect to the empirical $L_2$-norm, which results in a polynomial factor loss in the convergence. 
Its proof in Appendix \ref{prf_LowerBoundForTheSquaredExponentialKernel} follows from a more general lower bound in Proposition \ref{prp_LowerBound}.

\begin{figure}[H]
    \centering
    \begin{subfigure}{0.42\textwidth}
        \centering
        \includegraphics[width=\linewidth]{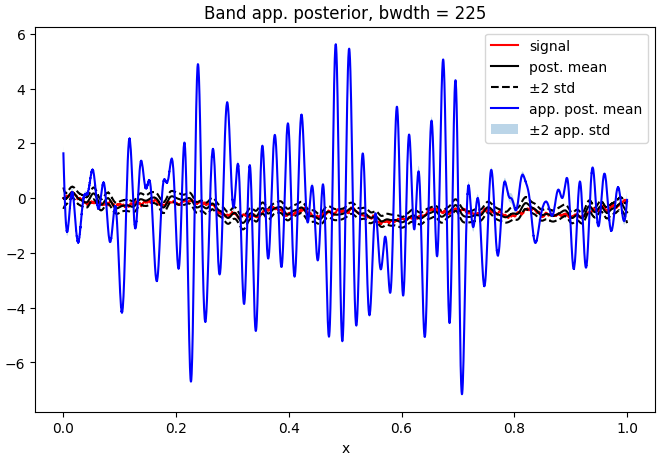}
        \caption{Approx. posterior with insufficient bandwidth.
                 Approx. posterior covariance is numerically not well defined.
                }
        \label{fig_squaredExponential_1}
    \end{subfigure}
    \hspace{0.5cm}
    \begin{subfigure}{0.42\textwidth}
        \centering
        \includegraphics[width=\linewidth]{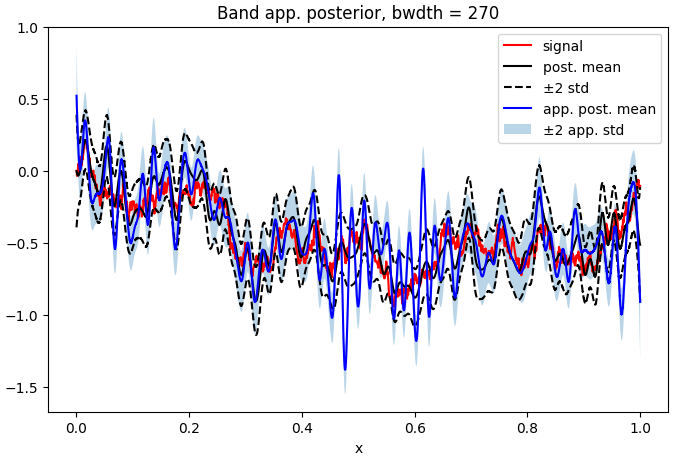}
        \caption{Approx. posterior with insufficient bandwidth.
                 Approx. posterior is well defined but deviates from the true one.
                }
        \label{fig_squaredExponential_2}
    \end{subfigure}

    \vspace{0.5cm}

    \begin{subfigure}{0.42\textwidth}
        \centering
        \includegraphics[width=\linewidth]{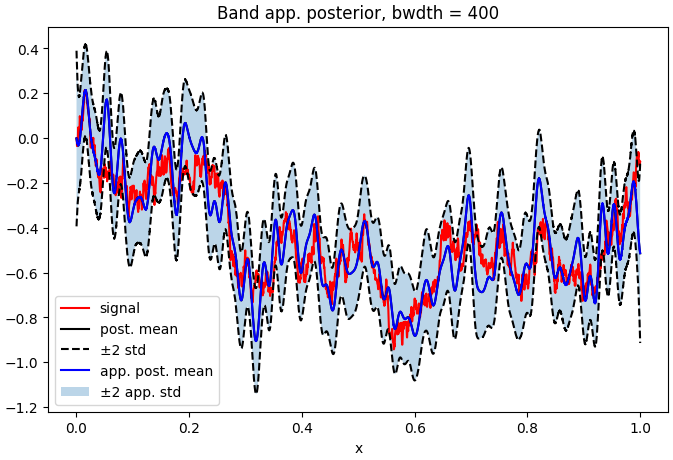}
        \caption{Approx. posterior with good bandwidth.
                }
        \label{fig_squaredExponential_3}
    \end{subfigure}
    \hspace{0.5cm}
    \begin{subfigure}{0.42\textwidth}
        \centering
        \includegraphics[width=\linewidth]{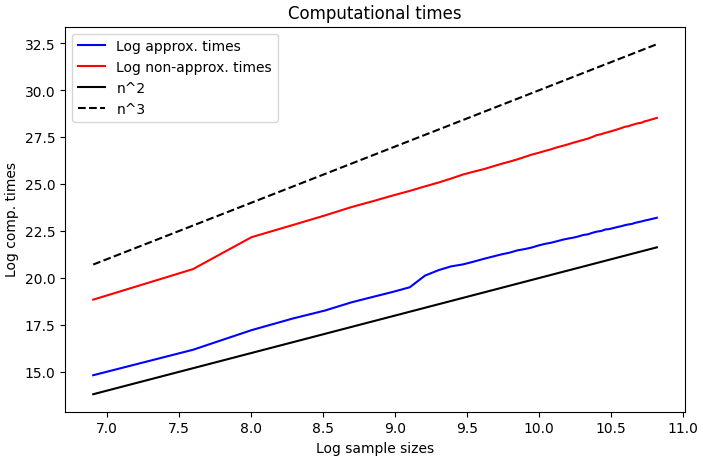}
        \caption{Log-log plot of the computation times.}
        \label{fig_squaredExponential_4}
    \end{subfigure}
    \caption{Squared exponential kernel simulations.
             The true function $f_0$ was set to a draw from a Brownian path and the covariates are iid uniformly distributed.
             In the simulations $n=10^4$, $\sigma=1$, and $\tau_n=n^{-1/(2\alpha+d)}$ were taken and the bandwidth $m$ was set to $225, 270$ and \( 400 \), respectively.
             }
    \label{fig_squaredExponential}
\end{figure}

Next, we demonstrate the performance of the approximation method numerically on synthetic data. 
We take a draw from a Brownian motion path as the true function $f_0$ and set $\sigma=1$ in the regression model from Equation \eqref{eq_Intro_GPRegressionModel}.  
For simplicity, we illustrate our simulations based on the choice of the number of neigbours $m$ instead of the truncation radius \( r \) noting that essentially \( m \asymp n r^d \).
We generate $n=10^4$ samples and consider the  number neighbours $m$ ranging from 225 to $400$, see Figure \ref{fig_squaredExponential}. 
Note that for a too small choice of $m$ the approximation is very poor, the process might not even be well defined (the corresponding covariance kernel is not positive definite), see Figure \ref{fig_squaredExponential_1}. 
As $m$ increases the resulting approximation becomes a valid Gaussian Process, but the approximation accuracy might not be acceptable, see Figure \ref{fig_squaredExponential_2}. 
Then after reaching a threshold the approximating GP, see Figure \ref{fig_squaredExponential_3} becomes almost indistinguishable from the true posterior and further increasing $m$ does not improve the approximation but only adds to the computational complexity. 
The log-computation time of the true posterior and the approximation is reported in Figure \ref{fig_squaredExponential_4}. 
The results are further illustrated in Table \ref{table:SE}, which reports the average and standard deviation of the squared $L_2$-error over 100 draws for different bandwidths $m$.

\begin{table}[H]
  \centering
  \begin{tabular}{ccccc}
    \\[-1.8ex]\hline
    \hline\\[-1.8ex]
    Bandwidth & 225 & 270 & 400 & \( 10^4 \) (full posterior) \\
    \hline\\[-1.8ex]
    Avg. MSE  & 159 & 32  & 0.01  & 0.01 \\
    Std. Dev. & 813 & 230 & 0.001 & 0.001 \\
    \\[-1.8ex]\hline
    \hline\\[-1.8ex]
  \end{tabular}
  \caption{Squared error of the true and approximate posterior based on SE kernel. Average and standard deviation over 100 repetition of the experiment.}\label{table:SE}
\end{table}


\subsection{The Mat\'ern kernel}
\label{ssec_TheMaternKernel}

In this  subsection, we discuss the Mat\'ern kernel. 
The Mat\'ern process is the centered GP with covariance kernel
\begin{align}
  \label{eq_Examples_MaternKernel}
  k_{\nu, \tau}(x, x') = \frac{2^{1-\nu}}{\Gamma(\nu)}\Big(\frac{\sqrt{2\nu}}{\tau} \| x' - x \|_2\Big)^\nu B_\nu\Big(\frac{\sqrt{2\nu}}{\tau} \| x' - x \|_2\Big), \quad x, x' \in\mathbb{R}^d,
\end{align}
where $B_\nu$ denotes the modified Bessel function of the second kind and $\Gamma(\nu)$ is the Gamma function.
When $\nu$ is half-integer, i.e. $\nu = \ell + 1/2$, where $\ell \in \mathbb N_0$, the Mat\'ern kernel has an explicit expression 
\begin{equation}
  \label{eq_Examples_MaternKernelHalfInteger}
  k_{\nu, \tau}(x, x')
  = 
  \exp\left(\frac{-\sqrt{2\ell+1} \| x' - x \|_2}{ \tau} \right) \frac{\Gamma(\ell + 1)}{\Gamma(2\ell+1)} \sum_{q=0}^\ell \frac{(\ell+q)!}{q!(\ell-q)!}\left( \frac{\sqrt{8\ell+4} \| x' - x \|_2}{\tau} \right)^{\ell - q},
\end{equation}
see Equation (4.16) in \cite{RasmussenEtal2006_GaussianProcesses}.
The sample paths from the process are $\nu$-Hölder smooth. A special case of the Mat\'ern process with $\nu=1/2$ is the Ornstein-Uhlenbeck process, also referred to as the GP with exponential kernel.

In the regression model one can recover the unknown $\alpha$-smooth function $f_0$ optimally by choosing a Mat\'ern process with parameter $\nu = \alpha$. 
In case the regularity of the process is not matching the regularity of $f_0$ one has to rescale the process for optimal recovery. For $f_0 \in C^{\alpha}(\mathbb{R}^d)$ with $0 < \alpha \le \ell + 1/2$ the optimal rescaling is given by \( \tau_n = n^{- (2 \ell + 1 - 2 \alpha) / ((2 \alpha + d) (2 \ell + 1))} \), see for instance \cite{FangBhadra2025_PosteriorConcentrationGaussianProcess, SzaboEtal2025_AdaptationSpatiallyDistributedGaussianProcesses}. 
Below we consider the Mat\'ern process prior with the above optimal rescaling and derive guarantees for the proposed approximation procedure.

\begin{proposition}[Result for the Mat\'ern kernel]
  \label{prp_ResultForTheMaternKernel}
  In model \eqref{eq_Intro_GPRegressionModel}, assume that \( \mathcal{X} \subset \mathbb R^d \) is bounded and the design observations are distributed according to \( G \) on \( \mathcal{X} \) with bounded Lebesgue density \( g \).
  Let \( f_0 \in C^{\alpha}(\mathcal X)\cap H^{\alpha}(\mathcal X) \) for \( \alpha > d/2 \) and let \( k_{\nu, \tau_n} \) be the rescaled Mat\'ern kernel with smoothness \( \nu = \ell + 1/2 \) for \( \ell \in \mathbb N_0 \) and lengthscale \( \tau_n = n^{- (\nu - \alpha) / (2 \alpha \nu + d \nu)} \).
  Then, the approximate posterior with truncation radius \( r_n = 2 \tau_n \log n  \) satisfies that for arbitrary \( M_n \to \infty \), both
  \begin{align}
    \Psi_{r_n}(d_H(\cdot, f_0) \ge M_n n^{- \alpha / (2 \alpha + d)}) \to 0
    \quad \text{ and } \quad 
    \kl(\Psi_{r_n}, \Pi_n(\cdot | X, Y)) \to 0
  \end{align}
  in probability with respect to \( \mathbb P_{f_0}^{\otimes n} \) for \( n \to \infty \).
  Further, the number of non-zero entries per row of \( \widehat K \) is at most of order \( n r_n^d \asymp n^{(2 \alpha + \alpha d / \nu) / (2 \alpha + d)} \log^{d} n\) with probability converging to one.
\end{proposition}

\noindent Proposition \ref{prp_ResultForTheMaternKernel}, which is proven in Appendix \ref{prf_ResultForTheMaternKernel}, shows that if the truncation radius is chosen at least of order \( n^{- (\nu - \alpha) / (2 \alpha \nu + d \nu)} \), up to logarithmic factors, then the approximate posterior attains the optimal contraction rate, and the Kullback-Leibler divergence between the approximate and true posteriors converges to zero.
Similarly to the squared exponential kernel, there is no substantial gap between the minimal choice of $ r_n $ required for vanishing approximation of the posterior measure in Kullback-Leibler divergence and obtaining the minimax optimal contraction rate.
This is again due to the exponential decay of the Mat\'ern kernel in the spatial domain. 
It is also worth emphasizing that the exponent \( -(\nu - \alpha) / (2 \alpha \nu + d \nu) \) is decreasing in \( \nu \) and increasing in \( \alpha \). 
Consequently, smoother kernels (larger \( \nu \)) allow for a smaller truncation radius \( r_n \).
Hence, considering smoother processes and rescaling them would allow to imply faster approximations. 
In contrast, as the smoothness \( \alpha \) of the true function \( f_0 \) increases, a larger neighborhood size is required. 
Finally, as \( \nu \) increases, the Mat\'ern kernel approaches the squared exponential kernel, and the corresponding minimal choice of \( r_n \) becomes comparable to that obtained in the SE setting up to logarithmic factors.
The same is true for the row-wise number \( m_n \asymp n r_n^d \) of non-zero entries in \( \widehat K \).
Again, storage is reduced to \( O(m_n n) \) and computational cost to \( O(m_n n N_{\text{iter}}) \).

\begin{figure}[H]
    \centering
    \begin{subfigure}{0.42\textwidth}
        \centering
        \includegraphics[width=\linewidth]{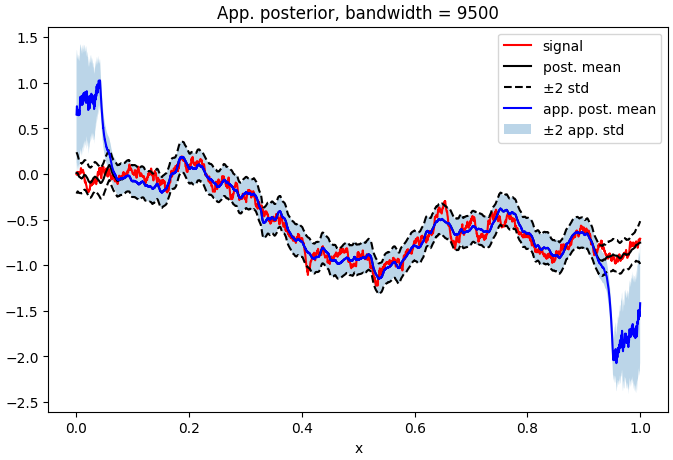}
        \caption{Regularity \( \nu = 0.5 \) with  \( \tau_{n} = 1 \).}
        \label{fig_MaternSimulations_1}
    \end{subfigure}
    \hspace{0.2cm}
    \begin{subfigure}{0.42\textwidth}
        \centering
        \includegraphics[width=\linewidth]{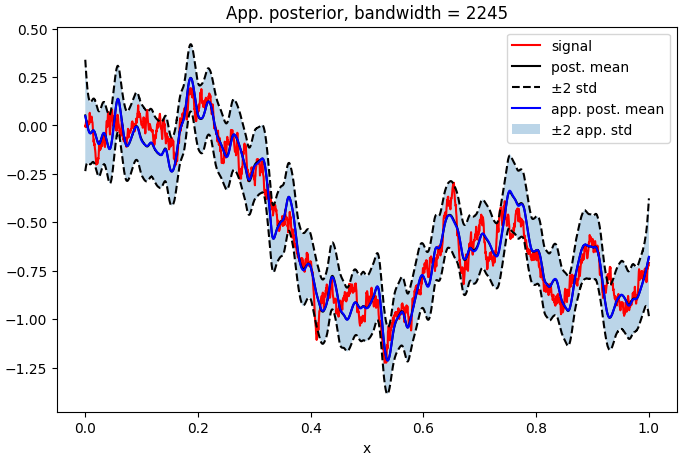}
        \caption{Regularity \( \nu = 1.5 \) with  \( \tau_{n} = n^{- 1/3} \).}
        \label{fig_MaternSimulations_2}
    \end{subfigure}

    \vspace{0.5cm}

    \begin{subfigure}{0.42\textwidth}
        \centering
        \includegraphics[width=\linewidth]{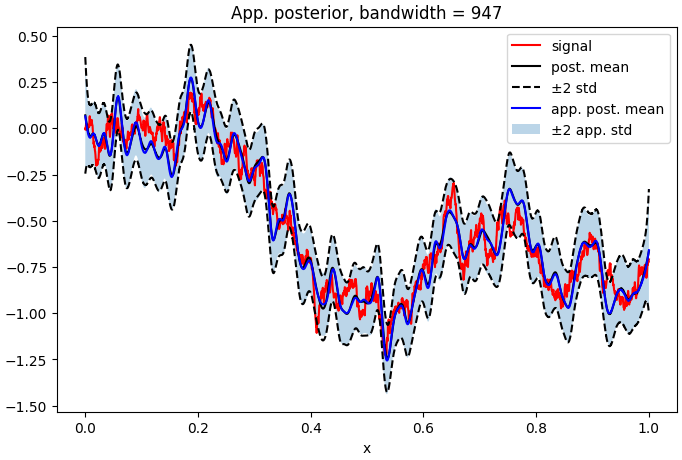}
        \caption{Regularity \( \nu = 2.5 \) with \( \tau_{n} = n^{- 2/5} \).}
        \label{fig_MaternSimulations_3}
    \end{subfigure}
    \hspace{0.2cm}
    \begin{subfigure}{0.42\textwidth}
        \centering
        \scriptsize
        \begin{tabular}{lccc}
          \\[-1.8ex]\hline
          \hline\\[-1.8ex]
          \( (\nu, \text{bwdth}) \) & \( (0.5, 9500) \) & \( (1.5, 2245) \) & \( (2.5, 947) \) \\
          \hline\\[-1.8ex]
          App. MSE             & 0.08              & 0.007             & 0.008            \\
          App. StD             & 0.5               & 0.001             & 0.001            \\
                               &                   &                   &                  \\
          True MSE             & 0.005             & 0.007             & 0.007            \\
          True StD             & 0.0007            & 0.001             & 0.001            \\
          \\[-1.8ex]\hline
          \hline\\[-1.8ex]
        \end{tabular}
        \vspace{20pt}
        \caption{Average mean squared error of the approximate and true posteriors with their standard deviations over 100 draws.}
        \label{fig_MaternSimulations_4}
    \end{subfigure}
    \caption{Simulations for the Matérn kernel.
             The true function $f_0$ was set to a draw from a Brownian motion path.
             In the simulation $n=10^4$ and $\sigma=1$ was taken and the number of neighbours $m$ was set to $9500, 2245 $ and \( 947 \), respectively.
            }
    \label{fig_MaternSimulations}
\end{figure}

Next we investigate the numerical performance of the approximation over a synthetic data set.
We study the influence of the regularity of the Mat\'ern process and the optimal choice of the neighbourhood size used in the approximation.
As in Section \ref{ssec_TheSquaredExponentialKernel}, we fix \( m \) explicitly instead of \( r \).
We set the true function $f_0$ to be a random draw of a Brownian path, take sample size $n=10^4$ and set $\sigma=1$ in the regression model from Equation \eqref{eq_Intro_GPRegressionModel}.
In Figure \ref{fig_MaternSimulations},  we consider the rescaled Mat\'ern process with regularities $\nu=1/2$, $\nu=3/2$ and $\nu=5/2$.
For \( \nu = 0.5 \), the optimal rescaling  is \( \tau_{n} = 1 \).
Consequently, according to our theory structural computational gains via banding should not be possible.
This is illustrated by the fact that even when keeping \( m = 9500 \) neighbours, the approximate and full posteriors  are still substantially different.
For \( \nu = 1.5 \) and the corresponding rescaling \( \tau_{n} = n^{- 1/3} \), we take \( m = 2245 \), the smallest neighbourhood size after which the approximation didn't improve further in our simulations (and became indistinguishable from the true posterior).
Finally, for \( \nu = 2.5 \) with rescaling \( \tau_{n} = n^{-2/5} \) this effect intensifies and the true and approximate posteriors already visually coincide already at \( m = 947 \).
The results are further confirmed by comparing the MSEs in Table (d).



\section{Discussion}
\label{sec_Discussion}

The results for both the squared exponential and the Matérn kernels exhibit a notable alignment between statistical optimality and approximation accuracy.
The order of the minimal truncation radius \( r = r_n \) required to achieve the optimal contraction rate and the one ensuring vanishing Kullback-Leibler divergence do not seem to differ substantially.
However, this phenomenon is not generic. 
As shown in Appendix \ref{ssec_APolynomialKernel}, for kernels with polynomial decay, the two requirements can differ at a polynomial scale in \(n\), showing that optimal recovery of the underlying function can be achieved with lower computational complexity than the one required for accurate approximation of the posterior.

We note that as the dimension $d$ increases, the required number of neighbours \( m_n \) in the approximation decreases in the above examples (for the Matérn case, this is true as long as \( \nu > \alpha \)).
In contrast, many other approximation methods exhibit complementary behaviour: as $d$ increases the scalability of the procedures reduces. 
For instance, in the inducing variable variational Bayes method, one has to involve $n^{d / (d + 2\alpha)}$ inducing variables \cite{NiemanEtal2022_ContractionRates, NiemanEtal2023_UncertaintyQuantificationSparseVariational}, the conjugate gradient descent and Lanczos approximations need to run for $O(n^{d / (d + 2 \alpha)})$ iterations \cite{StankewitzSzabo2024_ContractionRatesConjugateGradient}, and in the Vecchia approximation one requires parent sets of size ${\lfloor \alpha \rfloor + d \choose d}$ see \cite{SzaboZhu2026_VecchiaGaussianProcesses}. 
At the same time as the regularity $\alpha$ of the true function $f_0$ increases, one needs to keep more neighbours hence resulting in a less scalable procedure.

It is worth noting that the proposed banding strategy fundamentally relies on the kernel exhibiting a notion of spatial decay, typically of the form \( k(x, x') = r(\|x' - x \|) \), where \( r \) denotes a radial function. 
In this case, Euclidean proximity is aligned with correlation decay, which justifies the nearest-neighbour truncation. 
For general non-stationary kernels that do not admit such a representation, such as Brownian motion kernel, spatial distance  may therefore be inappropriate. 

Finally, the optimal choice of the truncation radius \( r \) depends on the unknown regularity \( \alpha \) of the true function. 
Since this quantity is not available in practice, an important direction for future work is to develop adaptive procedures that select \( r \) in a data-driven manner, for instance using early stopping. 



\begin{appendices}

\section{Proofs of the main results}
\label{sec_ProofsOfTheMainResults}

\begin{proof}[\normalfont \textbf{Proof of Theorem} \ref{thm_KullbackLeiblerBound}]
  \label{prf_KullbackLeiblerBound}
  The proof is based on the technical bound in Proposition \ref{prp_TechnicalKullbackLeiblerBound}.
  It suffices to bound the terms on the right-hand side with probability converging to one under \( \mathbb P_{f_0}^{\otimes n} \) as \( n \to \infty \).
  From Markov's inequality, in this sense, we have \( \| f_0 - h \|_n^2 \le M_n \| f_0 - h \|_{L^2}^2 \) for any sequence \( M_n \to \infty \).
  Now, under the concentration function assumption in Equation \eqref{eq_Main_PerturbationCondition}, \( h \) can be choseen such that \( \| h \|_{\mathbb H}^2 \lesssim n \eta_n^2 \) and \( \| f_0 - h \|_{L^2}^2 \lesssim \eta_n^2 \).
  From there, we obtain
  \begin{align}
            \frac{  \| h \|_\mathbb{H}^2 \| E \|_{\text{op}}^2} {\sigma^2 (\sigma^{2} - \| E \|_{\text{op}}) } 
          + 
            \frac{n \| f_0 - h \|_n^2    \| E \|_{\text{op}}^2} {\sigma^4 (\sigma^{2} - \| E \|_{\text{op}}) }
    & \le 
          \mathfrak C M_n \| E \|_{\text{op}}^2 n \eta_{n}^2.
  \end{align}
  Moreover, 
  \begin{align}
          \sum_{j = 1}^n \Big( \frac{\widehat \mu_j}{\widehat \mu_j + \sigma^2} \Big)^2  
    & \le 
          \sum_{j = 1}^n       \frac{\widehat \mu_j}{\widehat \mu_j + \sigma^2}
    = 
          \tr K K_\sigma^{-1}, 
  \end{align}
  which finishes the proof.
\end{proof}

\begin{proposition}[Technical Kullback-Leibler bound]
  \label{prp_TechnicalKullbackLeiblerBound}
  For any \( h \in \mathbb{H} \) and any sequence \( M_{n} \to \infty \), the approximate posterior satisfies\footnote{2026-09-19-BS-TODO: Formualte conditions if any.}
  \begin{align}
        \kl(\Psi_{r}, \Pi(\cdot | X, Y)) \mathbf{1} & \{ \| E \|_{\text{op}} < \sigma^2 / 2 \}
    \le 
             \frac{2   \| h \|_\mathbb{H}^2 \| E \|_{\text{op}}^2} {\sigma^2 (\sigma^{2} - \| E \|_{\text{op}}) } 
           + 
             \frac{2 n \| f_0 - h \|_n^2      \| E \|_{\text{op}}^2} {\sigma^4 (\sigma^{2} - \| E \|_{\text{op}}) }
    \\
    &   + 
             \frac{4 M_{n} \| E \|_{\text{op}}^2}{\sigma^4}
             \sum_{j = 1}^n \frac{\widehat \mu_j}{\widehat \mu_j + \sigma^2} 
           + 
             \frac{2 \| E \|_{\text{op}}^2} { (\sigma^{2} - 2 \| E \|_{\text{op}})^{2} } 
             \sum_{j = 1}^{n} \Big( \frac{\widehat{\mu}_{j}}{\widehat{\mu}_{j} + \sigma^{2}} \Big)^{2}
    \notag
  \end{align}
  with probability converging to one under \( \mathbb{P}_{f_0}^{\otimes n} \) and \( n \to \infty \).
\end{proposition}

\begin{proof}[Proof]
  In view of Lemma \ref{lem_ProcessKLVsFiniteDimensionalKL}, \( \kl(\Psi_{n,m}, \Pi(\cdot | X, Y)) = \kl(N_1, N_2) \), where  $ N_i \sim \mathcal N(m_i, \Sigma_i) $ with $ m_1 = K \widehat{K}_{\sigma}^{-1} Y $, $ m_2 = K {K}_{\sigma}^{-1} Y $, $ \Sigma_1 = (K - K \widehat{K}_{\sigma}^{-1} K) $ and $ \Sigma_2 = (K - K K_{\sigma}^{-1} K) $.
  \( N_1 \) and \( N_2 \) are always absolutely continuous with respect to each other.
  Indeed, when \( K \) is non-singular, both \( \Sigma_1 \) and \( \Sigma_2 \) are non-singular as well.
  Otherwise, a short calculation shows that \( \ran \Sigma_1 = \ran \Sigma_2 = \ran K \) and absolute continuity follows from the fact that borth \( m_1 \) and \( m_2 \) are also elements of \( \ran K \).
  Therefore, the Kullback-Leibler divergence between the finite dimensional Gaussians can be expressed as
  \begin{align}
    \label{eq_A_KLDecomposition}
        \kl(N_1, N_2)
    & = 
        \frac{1}{2} \big( \underbrace{ \tr \Sigma_{2}^{-1} \Sigma_{1} - n                     }_{(\text{I})  }
                        + \underbrace{ (m_{2} - m_{1})^{\top} \Sigma_{2}^{-1} (m_{2} - m_{1}) }_{(\text{II}) }
                        + \underbrace{ \log \det \Sigma_1^{-1} \Sigma_{2}                    }_{(\text{III})}
                    \big),
  \end{align}
  where ``\( (\cdot)^{-1} \)'' denotes the generalized inverse when \( K \) is singular.
  We deal with the three terms on the right hand side separately and restrict to the event \( \{ \| E \|_{\text{op}} < \sigma^2 / 2 \} \).

  \textit{Step 1: bound for the mean term (II).}
  By substitution and recalling that $ Y = \mathbf f_0 +  \varepsilon $, with  $\varepsilon \sim \mathcal{N}(0,\sigma^2 I_n)$ and $\mathbf f_0 = (f(X_1), \dots, f(X_n))^{\top} \in \mathbb{R}^n $, and the inequality $ \| a + b \|_2^2 \leq 2 \|a\|_2^2 + 2 \|b\|_2^2 $,
	\begin{align}
    \label{eq_ProofsMain_BoundMeanTermII}
           (\text{II}) 
    & =
           (K (K_{\sigma}^{-1} - \widehat{K}_{\sigma}^{-1}) Y)^{\top} (K - K K_{\sigma}^{-1} K)^{-1} 
           (K (K_{\sigma}^{-1} - \widehat{K}_{\sigma}^{-1}) Y) 
    \\
    & \le 
         2 (K (K_{\sigma}^{-1} - \widehat{K}_{\sigma}^{-1}) \varepsilon )^{\top} (K - K K_{\sigma}^{-1} K)^{-1} 
           (K (K_{\sigma}^{-1} - \widehat{K}_{\sigma}^{-1}) \varepsilon) 
    \notag
    \\
    & + 
           2 (K (K_{\sigma}^{-1} - \widehat{K}_{\sigma}^{-1}) \mathbf f_0)^{\top} (K - K K_{\sigma}^{-1} K)^{-1} 
             (K (K_{\sigma}^{-1} - \widehat{K}_{\sigma}^{-1}) \mathbf f_0).
    \notag
	\end{align}
	Now, we bound each term separately.
  In view of Markov's inequality for arbitrary $M_n \rightarrow \infty$ it holds with probability tending to one that
	\begin{align}
    & \ \ \ \
          \varepsilon^{\top} (K_{\sigma}^{-1} - \widehat{K}_{\sigma}^{-1}) 
          K (K - K K_{\sigma}^{-1} K)^{-1} 
          K (K_{\sigma}^{-1} - \widehat{K}_{\sigma}^{-1}) \varepsilon
    \\
		& \le 
          M_n \mathbb E_{\varepsilon} \tr \varepsilon \varepsilon^{\top} 
                                          (K_{\sigma}^{-1} - \widehat{K}_{\sigma}^{-1}) K
                                          (K - K K_{\sigma}^{-1} K)^{-1} 
                                          K (K_{\sigma}^{-1} - \widehat{K}_{\sigma}^{-1})
    \notag 
    \\
		& = 
          M_n \sigma^{2} \tr (K_{\sigma}^{-1} - \widehat{K}_{\sigma}^{-1}) K
                             (K - K K_{\sigma}^{-1} K)^{-1} 
                             K (K_{\sigma}^{-1} - \widehat{K}_{\sigma}^{-1})
    \notag 
    \\
    & =
          M_n \sigma^2 \tr \widehat K_{\sigma}^{-1} E K_{\sigma}^{-1} K
                           (K - K K_{\sigma}^{-1} K)^{-1} 
                           K K_{\sigma}^{-1} E\widehat K_{\sigma}^{-1} 
    \notag 
    \\
    & =
          M_n \sigma^2 \| (K - K K_{\sigma}^{-1} K)^{-1/2} K K_{\sigma}^{-1} E\widehat K_{\sigma}^{-1} \|_F^2
    \notag 
    \\
    & \le 
          M_n \sigma^2 \| (K - K K_{\sigma}^{-1} K)^{-1/2} K K_{\sigma}^{-1} \|_F^2 
                       \| E \|_{\rm op}^2 \| \widehat K_{\sigma}^{-1} \|_{\rm op}^2
    \notag 
    \\
    & \le 
          \frac{4 M_n \| E \|_{\rm op}^2}{\sigma^2} 
          \sum_{j = 1}^n \Bigl\{ \frac{\hat \mu_j}{\hat \mu_j+ \sigma^2 } \times
                                 \Bigl(\hat \mu_j- \frac{\hat \mu_j^2}{\hat \mu_j+ \sigma^2 } \Bigr)^{-1/2} 
                         \Bigr\}^2
    \notag 
    \\
    & = 
          \frac{4 M_n \| E \|_{\rm op}^2}{\sigma^4} 
          \sum_{j = 1}^n \frac{\hat \mu_j}{\hat \mu_j+ \sigma^2 },
    \notag 
	\end{align}
  where in the last inequality, we have used the fact that \( \| \widehat K_ \sigma^{-1} \|_{\text{op}} \le \sigma^2 / 2 \) when \( \| E \| < \sigma^2 / 2 \).

  For the second term in the bound \eqref{eq_ProofsMain_BoundMeanTermII},  we  split $ \mathbf f_0 =  (\mathbf f_0 - \mathbf h) + \mathbf h $, where for a function \( h \in \mathbb H \), $ \mathbf h = (h(X_1), \dots, h(X_n))^\top $. 
  Then, noting that
  \begin{align}
    K - K K_\sigma^{-1} K = K (I - K_\sigma^{-1} K) = \sigma^2 K K_\sigma^{-1}
  \end{align}
  and consequently
  \begin{align}
    K (K - K K_\sigma^{-1} K)^{-1} K = \sigma^{-2} K (K K_\sigma^{-1})^{-1} K
                                     = \sigma^{-2} K K_\sigma,
  \end{align}
  we obtain
	\begin{align}
    & \ \ \ \
            \mathbf h^{\top} (K_{\sigma}^{-1} - \widehat{K}_{\sigma}^{-1}) K 
                             (K - K K_{\sigma}^{-1} K)^{-1} 
                           K (K_{\sigma}^{-1} - \widehat{K}_{\sigma}^{-1})
            \mathbf h
    \\
    & = 
            \sigma^{-2} 
            \mathbf h^{\top} (K_{\sigma}^{-1} - \widehat{K}_{\sigma}^{-1}) 
                             K K_\sigma  
                             (K_{\sigma}^{-1} - \widehat{K}_{\sigma}^{-1})
            \mathbf h 
    \le 
            \sigma^{-2} 
            \| K_\sigma (K_{\sigma}^{-1} - \widehat{K}_{\sigma}^{-1}) \mathbf h \|^2
    \notag
    \\
    & = 
            \sigma^{-2} \| (I - K_\sigma \widehat{K}_{\sigma}^{-1}) \mathbf h \|^2
    = 
            \sigma^{-2} \| E \widehat{K}_{\sigma}^{-1} \mathbf h \|^2
    \le 
            \frac{\| E \|_{\text{op}}^2}{\sigma^2} \| \widehat{K}_{\sigma}^{-1} \mathbf h \|^2.
    \notag
	\end{align}
  From \( \widehat{K}_\sigma = K_\sigma + E \) it follows that \( \widehat{K}_\sigma^{-1} \le (1 - \| E \| / \sigma^2)^{-1} K_\sigma^{-1} \) and therefore
  \begin{align}
    & \ \ \ \
          \mathbf h^{\top} (K_{\sigma}^{-1} - \widehat{K}_{\sigma}^{-1}) K 
                           (K - K K_{\sigma}^{-1} K)^{-1} 
                         K (K_{\sigma}^{-1} - \widehat{K}_{\sigma}^{-1})
          \mathbf h
    \\
    & \le
          \frac{\| E \|_{\text{op}}^2}{ \sigma^2 - \| E \|_{\text{op}} } \| K_{\sigma}^{-1} \mathbf h \|^2
    \le 
          \frac{\| E \|_{\text{op}}^2}{\sigma^2 (\sigma^2 - \| E \|_{\text{op}})} 
          \mathbf h^\top K_{\sigma}^{-1} \mathbf h
    \le 
          \frac{\| E \|_{\text{op}}^2 \| h \|_{\mathbb H}^2}{\sigma^2 (\sigma^2 - \| E \|_{\text{op}})},
    \notag 
  \end{align}
  where we have used that \(\mathbf h^\top K_\sigma^{-1} \mathbf h \le \| h \|_{\mathbb H}^2 \) in the last inequality. 
  If \( K \) is non-singular, this follows from the fact that \( \mathbf h^{\top} K^{-1} \mathbf h \) is the squared RKHS norm of the projection of \( h \) onto \( \vspan \{ k(X_i, \cdot): i \le n \} \).
  Otherwise a simple additional calculation provides the same result.

  Dealing with the term $ \mathbf f_0 - \mathbf h$, by similar arguments,
	\begin{align}
    & \ \ \ \
          (\mathbf f_0 - \mathbf h)^\top (K_\sigma^{-1} - \widehat{K}_\sigma^{-1}) K 
                                       (K - K K_{\sigma}^{-1} K)^{-1} 
                                     K (K_\sigma^{-1} - \widehat{K}_\sigma^{-1}) 
          (\mathbf f_0 - \mathbf h)
    \\
    & \le 
          \frac{n \| f_0 - h \|_n^2 \| E \|_{\rm op}^2}{\sigma^4 (\sigma^2 -\| E \|_{\rm op})},
    \notag
	\end{align}
  where for a function \( g \), $\| g \|_n^2 = n^{-1} \sum_{i = 1}^n g(X_n)^2$.
  By combining the above displays, 
	\begin{equation}
        (\text{II}) 
    \le
        \frac{2 M_n \| E \|_{\rm op}^2}{\sigma^4} \sum_{j = 1}^n \frac{\hat \mu_j}{\hat \mu_j + \sigma^2}
      + 
        \frac{4   \| h \|_{\mathbb H}^2 \| E \|_{\rm op}^2}{\sigma^2 (\sigma^2 - \| E \|_{\rm op})}
      +
        \frac{4 n \| f_0 - h \|_n^2     \| E \|_{\rm op}^2}{\sigma^4 (\sigma^2 - \| E \|_{\rm op})}.
	\end{equation}

	\textit{Step 2: Decomposition of the log-determinant term (III).}
  Note that 
	\begin{align}
    & \ \ \ \
            (K - K K_{\sigma}^{-1} K)^{-1} (K - K \widehat K_\sigma^{-1} K) \\
    & = I + (K - K K_{\sigma}^{-1} K)^{-1} K (K_\sigma^{-1} - \widehat{K}_{\sigma}^{-1}) K
    \notag 
    \\
    & = I + (K - K K_{\sigma}^{-1} K)^{-1} K K_{\sigma}^{-1}E \widehat K_{\sigma}^{-1} K
      =  I + A, 
    \notag 
	\end{align}
  for \( A = (K - K K_{\sigma}^{-1} K)^{-1} K K_{\sigma}^{-1} E \widehat K_{\sigma}^{-1} K \).
  By Lemma \ref{lem_NormEstimateForTheLogDeterminantTerm} and the fact that \( \| E \|_{\rm op} < \sigma^2 / 2 \), we have \( \| A \|_{\rm op} \le \| E \|_{\rm op} / (\sigma^2 - \| E \|_{\rm op}) < 1 \).
	Thus, the matrix logarithm $\log(I_n + A)$ is well-defined and we obtain
  \begin{align}
    (\text{III}) = - \log \det (I_n + A) = - \tr \log (I_n + A)
  \end{align}

  The matrix logarithm admits the convergent (in operator norm) series
  \begin{align}
    \log(I_n + A) = \sum_{k = 1}^\infty \frac{(-1)^{k + 1}}{k} A^k.
  \end{align}
	Taking traces term by term yields
  \begin{align}
    \label{eq_ProofsMain_TraceDecompositionOfTheLogDeterminant}
    (\text{III}) = - \tr \log(I_n + A) 
                 = - \tr A + \frac{1}{2} \tr A^2 
                 - 
                   \sum_{k = 3}^\infty \frac{(-1)^{k + 1}}{k} \tr A^k.
  \end{align}
  In the same notation, the first term \( (\text{I}) \) can be expressed as
  \begin{align}
    (\text{I}) & = \tr (K - K K_{\sigma}^{-1} K)^{-1} (K - K \widehat{K}_{\sigma}^{-1} K) - n
                 = \tr A
  \end{align}
  such that
  \begin{align}
    (\text{I}) + (\text{III}) = \frac{1}{2} \tr A^2 
                              - \sum_{k = 3}^\infty \frac{(-1)^{k + 1}}{k} \tr A^k.
  \end{align}

	\textit{Step 3: The square trace term.}
  For the square trace term, we obtain
  \begin{align}
           | \tr A^2| 
    & = 
           | \tr   (K - K K_\sigma^{-1} K)^{-1} K K_\sigma^{-1} E \widehat K_\sigma^{-1} K 
             \cdot (K - K K_\sigma^{-1} K)^{-1} K K_\sigma^{-1} E \widehat K_\sigma^{-1} K
           |
    \\
    & = 
           | \tr   \widehat K_\sigma^{-1} K (K - K K_\sigma^{-1} K)^{-1} K K_{\sigma}^{-1} E  
             \cdot \widehat K_\sigma^{-1} K (K - K K_\sigma^{-1} K)^{-1} K K_{\sigma}^{-1} E
           |
    \notag
    \\
    & \le 
           \| \widehat K_\sigma^{-1} K (K - K K_\sigma^{-1} K)^{-1} K K_\sigma^{-1} E  \|_{\text{F}}^2
    \notag
    \\
    & \le 
           \| \widehat K_\sigma^{-1} K (K - K K_\sigma^{-1} K)^{-1} K K_\sigma^{-1} \|_{\text{F}}^2
           \| E \|_{\text{op}}^2,
    \notag
  \end{align}
  where we have used elementary estimates for the Frobenius norm and the fact that all matrices involved are symmetric.
  Next, we use the identity \( \widehat K_\sigma^{-1} = (K_\sigma + E)^{-1} = (I_n + K_ \sigma^{-1} E)^{-1} K_ \sigma^{-1} \), which yields the final estimate
  \begin{align}
           | \tr A^2| 
    & \le 
           \| \widehat K_\sigma^{-1} K (K - K K_\sigma^{-1} K)^{-1} K K_\sigma^{-1} \|_{\text{F}}^2
           \| E \|_{\text{op}}^2
    \\
    & \le  
           \| (I_n + K_ \sigma^{-1} E)^{-1} \|_{\text{op}}^2
           \| K_ \sigma^{-1} K (K - K K_\sigma^{-1} K)^{-1} K K_\sigma^{-1} \|_{\text{F}}^2 
           \| E \|_{\text{op}}^2
    \notag
    \\
    & = 
           \| (I_n + K_ \sigma^{-1} E)^{-1} \|_{\text{op}}^2
           \sum_{j = 1}^n \Big( \frac{\widehat \mu_j^2         }{(\widehat \mu_j + \sigma^2)^2} 
                                \frac{\widehat \mu_j + \sigma^2}{ \widehat \mu_j   \sigma^2   } 
                          \Big)^2
           \| E \|_{\text{op}}^2
    \notag
    \\
    & \le 
           \frac{\| E \|_{\text{op}}^2}{(\sigma^2 - \| E \|_{\text{op}})^2} 
           \sum_{j = 1}^n \Big( \frac{\widehat \mu_j}{\widehat \mu_j + \sigma^2} \Big)^2
    \notag
  \end{align}
  with the estimate \( \| (I_n + K_ \sigma^{-1} E)^{-1} \|_{\text{op}}^2 \le 1 / (1 - \| E \|_{\text{op}} / \sigma^2) \).

	\textit{Step 4: Controlling the remainder.}
  Let \( k \ge 3 \) be an integer.
  Then, we have by Cauchy-Schwartz inequality that
  \begin{align}
    & \ \ \ \
           |\tr A^{k}| 
    \\
    & = 
           \big| \tr \big( (K - K K_\sigma^{-1} K)^{-1} K K_\sigma^{-1} E \widehat K_\sigma^{-1} K \big)^k \big| 
    \notag
    \\
    & = 
           \big| \tr \big( \widehat K_\sigma^{-1} K (K - K K_\sigma^{-1} K)^{-1} K K_\sigma^{-1} E \big)^k \big| 
    \notag
    \\
    & \le 
           \big\| \big( \widehat K_\sigma^{-1} K (K - K K_\sigma^{-1} K)^{-1} K K_\sigma^{-1} E \big)^{k - 1} \big\|_{\text{F}} 
           \big\| \big( \widehat K_\sigma^{-1} K (K - K K_\sigma^{-1} K)^{-1} K K_\sigma^{-1} E \big)         \big\|_{\text{F}}
    \notag
    \\
    & \le 
           \big\| \big( \widehat K_\sigma^{-1} K (K - K K_\sigma^{-1} K)^{-1} K K_\sigma^{-1} E \big)^{k - 2} \big\|_{\text{op}} 
           \big\| \big( \widehat K_\sigma^{-1} K (K - K K_\sigma^{-1} K)^{-1} K K_\sigma^{-1} E \big)         \big\|_{\text{F}}^2
    \notag
  \end{align}
  As in Step 3, we have
  \begin{align}
        \big\| \big( \widehat K_\sigma^{-1} K (K - K K_\sigma^{-1} K)^{-1} K K_\sigma^{-1} E \big)         \big\|_{\text{F}}^2 
    \le 
        \frac{\| E \|_{\text{op}}^2}{(\sigma^2 - \| E \|_{\text{op}})^2}
        \sum_{j = 1}^n \Big( \frac{\widehat \mu_j}{\widehat \mu_j + \sigma^2} \Big)^2
  \end{align}
  Further, again using the identity \( \widehat K_\sigma^{-1} = (I_n + K_\sigma^{-1} E)^{-1} K_\sigma^{-1} \),  
  \begin{align}
    & \ \ \ \
            \big\| \big( \widehat K_\sigma^{-1} K (K - K K_\sigma^{-1} K)^{-1} K K_\sigma^{-1} E \big)^{k - 2} \big\|_{\text{op}} 
    \le 
            \big\|       \widehat K_\sigma^{-1} K (K - K K_\sigma^{-1} K)^{-1} K K_\sigma^{-1} E               \big\|_{\text{op}}^{k - 2} 
    \notag
    \\
    & \le 
                \| E \|_{\text{op}}^{k - 2} 
                \| (I_n + K_\sigma^{-1} E)^{-1} \|_{\text{op}}^{k - 2}
            \big\| K_\sigma^{-1} K (K - K K_\sigma^{-1} K)^{-1} K K_\sigma^{-1} \big\|_{\text{op}}^{k - 2} 
    \notag
    \\
    & \le 
            \frac{ \| E \|_{\text{op}}^{k - 2} }{ (1 - \| E \|_{\text{op}} / \sigma^2)^{k - 2} \sigma^{2 (k - 2)}} 
            \max_{j \le n} \Big( \frac{\widehat \mu_j}{\widehat \mu_j + \sigma^2} \Big)^{k - 2}
    \le 
            \Big( \frac{ \| E \|_{\text{op}} }{ (\sigma^2 - \| E \|_{\text{op}}) } \Big)^{k - 2}.
  \end{align}
  For the whole remainder in \( \text{(I)} + \text{(III)} \), we therefore obtain
  \begin{align}
    \sum_{k = 3}^\infty \frac{1}{k} | \tr A^k |
    & \le 
    \frac{\| E \|_{\text{op}}^2}{(\sigma^2 - \| E \|_{\text{op}})^2}
    \sum_{j = 1}^n \Big( \frac{\widehat \mu_j}{\widehat \mu_j + \sigma^2} \Big)^2
    \sum_{k = 3}^\infty 
    \frac{1}{k}
    \Big( \frac{ \| E \|_{\text{op}} }{ \sigma^2 - \| E \|_{\text{op}} } \Big)^{k - 2}
    \\
    & \le
    \frac{\| E \|_{\text{op}}}{(\sigma^2 - \| E \|_{\text{op}})}
    \sum_{j = 1}^n \Big( \frac{\widehat \mu_j}{\widehat \mu_j + \sigma^2} \Big)^2
    \sum_{k = 1}^\infty 
    \frac{1}{k}
    \Big( \frac{ \| E \|_{\text{op}} }{ \sigma^2 - \| E \|_{\text{op}} } \Big)^k
    \notag
    \\
    & = 
    \frac{- \| E \|_{\text{op}}}{(\sigma^2 - \| E \|_{\text{op}})}
    \sum_{j = 1}^n \Big( \frac{\widehat \mu_j}{\widehat \mu_j + \sigma^2} \Big)^2
    \log \Big( 1 - \frac{\| E \|_{\text{op}} }{ \sigma^2 - \| E \|_{\text{op}} } \Big)
    \notag
    \\
    & \le 
    \frac{\| E \|_{\text{op}}^2}{(\sigma^2 - 2 \| E \|_{\text{op}})^2}
    \sum_{j = 1}^n \Big( \frac{\widehat \mu_j}{\widehat \mu_j + \sigma^2} \Big)^2,
    \notag
  \end{align}
  where we have used the inequality \( - \log(1 - x) \le x / (1 - x) \) for \( x \in [0, 1) \) in the last step.
  Combining the results yields the final bound.
\end{proof}

\begin{proof}[\normalfont \textbf{Proof of Theorem \ref{thm_ContractionRates}}]
  \label{prf_ContractionRates}
  By Theorem \ref{thm_ContractionRates} and the additional assumption on \( \tr K K_ \sigma^{-1} \) there is an event with probability converging to one under \( \mathbb P_{f_0}^{\otimes n} \), such that for any \( M_n \to \infty \),
  \begin{align}
    \kl(\Psi_{n, m}(\cdot), \Pi_n(\cdot | \mathcal D^n)) \lesssim M_n n \eta_n^2.
  \end{align}
  The stattement of the theorem now follows by applying Proposition 6.1 in \cite{StankewitzSzabo2024_ContractionRatesConjugateGradient}, see also Theorem 5 in \cite{RaySzabo2022_VariationalBayesHighDimensionalRegression}.
\end{proof}


\section{Proofs for the Examples}
\label{sec_ProofsForTheExamples}

\begin{proof}[\normalfont \textbf{Proof of Proposition \ref{prp_ResultForTheSquardExponentialKernel}} (Result for the squared exponential kernel)]
  \label{prf_ResultForTheSquardExponentialKernel}
  Initially, we note that it is a well known fact that under the assumptions of Proposition \ref{prp_ResultForTheSquardExponentialKernel}, the concentration function inequality \eqref{eq_Main_ConcentrationFunctionInequality} is satisfied for \( \eta_n = n^{-\alpha / (2 \alpha + d)} \log^{(1/2 + 3d/4)} n \), see,  e.g., Lemma 8 in \cite{NiemanEtal2022_ContractionRates}.

  Consequently, in order to apply Theorems \ref{thm_ContractionRates} and \ref{thm_ContractionRates}, it suffices to bound \( \| E \|_{\text{op}} = \| \widehat K - K \|_{\text{op}} \) and \( \tr K K_\sigma^{-1} \) with probability converging to one.
  For \( r = r_n = \sqrt{2 \tau _n^2 \log n} \), we immediately have that
  \begin{align}
          \| \widehat K - K \|_{\text{op}} 
    & \le 
          \max_{i \le n} \sum_{j: \| X_j - X_i \| > r} k_{\tau_n}(X_i, X_j) 
    \le 
          n \exp \Big( \frac{-r_n^2}{\tau_n^2} \Big) 
    \le 
          \frac{1}{n} \to 0.
  \end{align}
  In particular, \( \| \widehat K - K \|_{\text{op}} < \sigma^2 / 2 \) eventually for all \( n \in \mathbb N \).
  Further, by Lemma \ref{lem_TraceTermForTheSquaredExponentialKernel}, for any \( M_n \to \infty \),
  \begin{align}
    \tr K K_ \sigma^{-1} \lesssim M_n n^{\frac{d}{2 \alpha + d}} \log^d n
                         \lesssim M_n n \eta_n^2
  \end{align}
  with probability converging to one under \( \mathbb P_{f_0}^{\otimes n} \) and \( n \to \infty \).
  Both the contraction rate and the convergence of the Kullback-Leibler divergence follow from combining these bounds.

  Finally, the statement about the non-zero row entries follows from applying Lemma \ref{lem_NumberOfNeighboursBound} with \( n r^d = n r_n^d \asymp n \tau_n^{d} \log^{d/2} n = n^{2 \alpha / (2 \alpha + d)} \log^{d/2} n \) and \( \delta = \| g \|_{\infty} r^d \text{Vol}(B_1(0)) \) together with a union bound.
  Note that \( \alpha > d/2 \) guarantees that \( n \delta^2 \to \infty \) in the exponential bound.
\end{proof}

\begin{proof}[\normalfont \textbf{Proof of Corollary \ref{cor_LowerBoundForTheSquaredExponentialKernel}} (Lower bound for the squared exponential kernel)]
  \label{prf_LowerBoundForTheSquaredExponentialKernel}
  We first consider the case $ d = 1$.
  Assume that $G = \mathbb P^X$ has density $g(x) \propto e^{-x^2 / 2}$.
  As in Lemma \ref{lem_TraceTermForTheSquaredExponentialKernel}, the eigenvalues of the covariance operator then satisfy
  \begin{align}
    \lambda_j = \sqrt{1 / (2 A_n)} \left(\frac{1}{A_n \tau_n^2}\right)^{j-1}, \qquad j \ge 1,
  \end{align}
  where $A_n = 1/4 + \tau_n^{-2} + \sqrt{1/16 + \tau_n^{-2} / 2}$.
  In particular \( \lambda_1 = \sqrt{1/(2 A_n)} \asymp \tau_n \).
  For $ d > 1$, by tensorisation \( \lambda_1 \asymp \tau_n^d \).
  Noting that the first eigenfunction of the kernel operator is given by \( \phi_1(x) = (1 + 4 a_\tau)^{d/4} \exp(-a_\tau \| x \|^2) \) with \( a_ \tau = (\sqrt{1 + 4/\tau^2} - 1) / (4 (1 + \tau^2)) \), see \cite{RasmussenEtal2006_GaussianProcesses}, we obtain that \( \varphi_1 = \lambda_{1}^{1/2} \phi_1 \in C^{\alpha}(\mathbb R^d) \) for any fixed \( \alpha > 0 \).
  Further, it is a well known fact that in this setting the squared exponential kernel is strictly positive definite in the sense of Proposition \ref{prp_LowerBound}.

  Next, for \( r_n = n^{-1 / (2 \alpha + d)} / \log^{1/(2d)} n \), by Lemma \ref{lem_NumberOfNeighboursBound} and the boundedness of the kernel,
  \begin{align}
             \|\widehat K\|_{\mathrm{op}} 
    \le      \max_{1 \le i \le n} \sum_{j: \| X_j - X_i \| \le r_n} k_{\tau_n}(X_i, X_j)
    \lesssim n r_n^d
  \end{align}
  with probability converging to one. 
  Taking $ m_n: = n r_n^d \asymp n^{2 \alpha / (2 \alpha + d)} / \sqrt{\log n} $ and using the fact that \( \tau_n = n^{-1 / (2 \alpha + d)} \), we obtain
  \begin{align}
            \inf_{r \le r_n} 
            \sqrt{\lambda_1} \left( \frac{n \lambda_1}{2\|\widehat K\|_{\mathrm{op}}} - 1 \right) 
    \gtrsim \sqrt{\tau_n^d} \left( \frac{n \tau_n^d}{m_n} \right) 
    \gtrsim n^{\frac{-d/2}{2 \alpha + d}} \sqrt{\log n}.
  \end{align}
  This controls the first tern in the bound from Proposition \ref{prp_LowerBound}. 

  We now control the second term using
  \begin{align}
    \mathbb{E}\|k(X_1,\cdot)\|_{\mathbb H}^4 \lesssim 1, 
    \qquad 
    \|\Sigma\|_{\mathrm{HS}}^2 = \sum_{j\ge1}\lambda_j^2  \asymp \tau_n^{d},
  \end{align}
  where the second statement follows from Equation \eqref{eq_Aux_EVDecayForSEKernel}.
  Hence, \( \mathbb{E} \| k(X_{1}, \cdot) \|_{\mathbb{H}}^{4} - \| \Sigma \|_{\text{HS}}^{2} \lesssim 1 \). 
  Moreover, using the exact formulas for $ \lambda_1 $ and $ \lambda_2 $, we obtain $\lambda_1 - \lambda_2 \asymp \tau_n^{d + 1} $.
  Therefore,
  \begin{align}
    \frac{\sqrt{\mathbb{E} \| k(X_{1}, \cdot) \|_{\mathbb{H}}^{4} - \| \Sigma \|_{\text{HS}}^{2}}}
         {(\lambda_{1} - \lambda_{2}) \sqrt{n / \log n}}
    \;\lesssim\;
    \frac{1}{\tau_n^{d + 1} \sqrt{n/\log n}}.
  \end{align}
  Multiplying by $\sqrt{\lambda_1} \asymp \tau_n^{d/2}$ and using $\tau_n = n^{-1/(2\alpha+d)}$, we get
  \begin{align}
    \sqrt{\lambda_1}
    \frac{\sqrt{\mathbb{E} \| k(X_{1}, \cdot) \|_{\mathbb{H}}^{4} - \| \Sigma \|_{\text{HS}}^{2}}}
         {(\lambda_{1} - \lambda_{2}) \sqrt{n / \log n}}
    \;\lesssim\;\sqrt{\log(n)}
    n^{\frac{- \alpha + 1}{2\alpha+d}}.
  \end{align}

  Combining both results and applying Proposition~\ref{prp_LowerBound}, we conclude that with probability tending to $1$, the approximate posterior is either not well defined or
  \begin{align}
    \inf_{r \le r_n} \sup_{f_0 \in C^{\alpha}(\mathbb R^d)} 
    \| \mathbb{E}_{\varepsilon} K \widehat{K}_{\sigma}^{-1} Y - f_{0} \|_{n}
    \gtrsim
    \Bigl( \sqrt{\log(n)}\, n^{\frac{-d}{2(2\alpha+d)}} -\sqrt{\log(n)}\, n^{\frac{- \alpha + 1}{2\alpha+d}} \Bigr)^{+} .
  \end{align}
  The result follows since $\alpha > d/2 + 1$.
\end{proof}

\begin{proof}[\normalfont \textbf{Proof of Proposition \ref{prp_ResultForTheMaternKernel}} (Result for the Matérn kernel)]
  \label{prf_ResultForTheMaternKernel}
  From Lemmas 3.2 and 3.3 in \cite{FangBhadra2025_PosteriorConcentrationGaussianProcess}, it follows that the concentration function inequality \eqref{eq_Main_ConcentrationFunctionInequality} is satisfied with \( \eta_n = n^{-\alpha / (2 \alpha + d)} \).  

  As for the squared exponential kernel, it is then sufficient to bound \( \| E \|_{\text{op}} = \| \widehat K - K \|_{\text{op}} \) and \( \tr K K_\sigma^{-1} \) with probability converging to one.
  We use the fact that for \( \nu = \ell + 1/2 \), we have an explicit representation of the Mat\'ern kernel such that with \( r = r_n = 2 \tau_n \log n \), we have
  \begin{align}
    \| \widehat K - K \|_{\text{op}} 
    & \le 
    \max_{i \le n} \sum_{j: \| X_j - X_i \| > r_n} k_{\tau_n}(X_i, X_j)
    \\
    & \lesssim
    \max_{i \le n} \sum_{j: \| X_j - X_i \| > r_n} 
                   \frac{\ell !}{(2 \ell)!} \exp \Big( \frac{- \sqrt{2 \ell + 1} r_n}{\tau_n} \Big) 
                   \sum_{q = 0}^{\ell}
                   \frac{(\ell + q)!}{q! (\ell - q)!} \Big( \frac{2 \sqrt{2 \ell + 1} r_n}{\tau_n} \Big)^{\ell - q}
    \notag
    \\
    & \le 
    n \exp \big( -2 \sqrt{2 \ell + 1 } \log n \big)
    \sum_{q = 0}^{\ell} \binom{\ell}{q} \big( 4 \sqrt{2 \ell + 1} \log n \big)^{\ell - q}
    \notag
    \\
    & =
    n \exp \big( -2 \sqrt{2 \ell + 1 } \log n \big)
    \big( 4 \sqrt{2 \ell + 1} \log n + 1 \big)^{\ell}
    \notag
    \\
    & \lesssim 
    n^{1 - \sqrt{8 \ell + 4}} \log^{\ell} n 
    \to 0
    \notag
  \end{align}
  where we have used the binomial theorem in the last equality.
  The required bound on \( \tr K K_\sigma^{-1} \) is provided in Lemma \ref{lem_TraceTermForTheMaternKernel}.

  Finally, the statement about the non-zero row entries follows from applying Lemma \ref{lem_NumberOfNeighboursBound} with \( n r^d = n r_n^d \asymp n \tau_n^{d} \log^d n = n^{(2 \alpha \nu + d \alpha) / (2 \alpha \nu + d \nu) } \log^d n \) and \( \delta = \| g \|_{\infty} r^d \text{Vol}(B_1(0)) \) together with a union bound.
  Note that \( \alpha > d/2 \) guarantees that \( n \delta^2 \to \infty \) in the exponential bound for any choice of \( \nu = l + 1/2 \) with \( \ell \in \mathbb N_0 \).
\end{proof}


\section{Auxiliary results}
\label{sec_AuxiliaryResults}

\begin{lemma}[Norm estimate for the log-determinant term]
  \label{lem_NormEstimateForTheLogDeterminantTerm}
  Set $ E = \widehat K - K $.
  Then 
  \begin{align}
        \| (K - K K_{\sigma}^{-1} K)^{-1} K K_{\sigma}^{-1} E \widehat K_{\sigma}^{-1} K \|_{\rm op} 
    \le
        \frac{ \| E \|_{\rm op} }{ \sigma^2 - \| E \|_{\rm op} }.
  \end{align}
\end{lemma}

\begin{proof}
  Since $\widehat K_{\sigma} = K_{\sigma}+E$, we have
  \begin{align}
    & \ \ \ \
    \|(K - K K_{\sigma}^{-1} K)^{-1} K K_{\sigma}^{-1} E \widehat K_{\sigma}^{-1} K \|_{\rm op} 
    \\
    & = 
    \| (K - K K_{\sigma}^{-1} K)^{-1} K K_{\sigma}^{-1} E (I_n + K_{\sigma}^{-1} E)^{-1} K_{\sigma}^{-1} K \|_{\rm op}
    \notag
    \\
    &\le
    \|(K-KK_{\sigma}^{-1}K)^{-1}K K_{\sigma}^{-1} \|_{\rm op} \,\|E\|_{\rm op} \,
    \|(I_n+K_{\sigma}^{-1}E)^{-1}\|_{\rm op} \,\|K_{\sigma}^{-1}K\|_{\rm op}.
    \notag 
  \end{align}
  We now bound each factor individually.
  First, \( \| K_{\sigma}^{-1} K \|_{\text{op}} \le 1 \).

  Next, since $ \| K_{\sigma}^{-1} E \|_{\rm op} < 1$, the Neumann series gives
  \begin{align}
    \| (I_n + K_{\sigma}^{-1} E)^{-1} \|_{\rm op} \le \frac{1}{ 1 - \|K_{\sigma}^{-1} E\|_{\rm op} }.
  \end{align}
  Moreover, \( \| K_{\sigma}^{-1} E \|_{\rm op} \le \| K_{\sigma}^{-1} \|_{\rm op} \| E \|_{\rm op} \le \| E \|_{\rm op} / \sigma^2 \).
  Hence,
  \begin{align}
    \| (I_n + K_{\sigma}^{-1} E)^{-1} \|_{\rm op} \le \frac{1}{1 - \| E \|_{\rm op} / \sigma^2}.
  \end{align}
  Finally, since $K$ and $K_{\sigma}$ are simultaneously diagonalizable, we obtain
  \begin{align}
        \| (K - K K_{\sigma}^{-1} K)^{-1} K K_{\sigma}^{-1} \|_{\rm op}
    & =
        \max_{j = 1, \dots, n} \Big\{ \Big( \hat \mu_j - \frac{\hat \mu_j^2}{\hat \mu_j +\sigma^2} \Big)^{-1}
                                      \frac{\hat \mu_j}{\hat \mu_j +\sigma^2}
                              \Big\}
      = \sigma^{-2}.
  \end{align}
  Combining the previous bounds yield the result.
\end{proof}

\begin{lemma}[Trace term for the squared exponential kernel]
  \label{lem_TraceTermForTheSquaredExponentialKernel}
  Assume that \( \mathcal{X} = \mathbb{R}^{d} \) and that \( G = N(0, \gamma^2 I_d) \) for some \( \gamma^2 > 0 \).  
  Let \( k_{\tau_n} \) be the squared exponential kernel with length scale \( \tau_n = n^{-1/(2 \alpha + d)} \), where \( \alpha > 0 \).  
  Then, for any sequence \( M_n \to \infty  \), 
  \begin{align}
    \tr(K K_\sigma^{-1}) \lesssim M_n n^{\frac{d}{2 \alpha + d}} \log^d n,
  \end{align}
  with probability converging to one under \( \mathbb{P}_{f_0}^{\otimes n} \) as \( n \to \infty \).
\end{lemma}

\begin{proof}[Proof]
	First we give an upper bound for the eigenvalues for the squared exponential covariance kernel.
  This derivation is standard, but for completeness we provide a proof.	
  For $ d=1 $ case, let Let $ (\phi_j)_{j \ge 1} $ and $(\lambda_{1, j})_{j \ge 1}$ denote the eigenfunctions and eigenvalues of the one-dimensional kernel $k_{\tau_n}$,  respectively.
  For notational convenience, assume that $ G $ has a density $g(x) \propto e^{-2 a x^2}$ for some $a>0$.
  In this case, the eigenvalues admit the explicit form
  \begin{align}
    \lambda_{1, j} = \sqrt{2a/A_n} \left( \frac{1}{A_n \tau_n^2}\right)^{j - 1}, \qquad j \ge 1,
  \end{align}
  where $A_n = a + \tau_n^{-2} + \sqrt{a^2 + 2 a \tau_n^{-2}}$, see \cite{RasmussenEtal2006_GaussianProcesses}.
  With elementary derivations, it can then be shown that \begin{align}
    \label{eq_Aux_EVDecayForTheSEKernelD1}
    \lambda_{1, j} \lesssim \tau_n \exp(-D_{\exp} \tau_n j), \qquad j \ge 1
  \end{align}
  for any \( D_{\exp} < \sqrt{2 a} \), see, e.g., the proof of Corollary 9 in \cite{NiemanEtal2022_ContractionRates}.

  Next, the eigenfunctions and eigenvalues of the general \( d \)-dimensional kernel \(k_{\tau_n} \) are given by
  \begin{align}
    \Phi_{\mathbf{j}}(x) = \prod_{k = 1}^d \phi_{j_k}(x_k),  \quad
    \Lambda_{\mathbf{j}} = \prod_{k = 1}^d \lambda_{1, j_k}, \qquad
                           \mathbf{j}=(j_1, \dots, j_d)\in \mathbb{N}^d.
  \end{align}
  Using the one-dimensional bound in Equation \eqref{eq_Aux_EVDecayForTheSEKernelD1}, we obtain that for every $\mathbf{j}\in\mathbb{N}^d$,
  \begin{align}
    \Lambda_{\mathbf{j}} =              \prod_{k = 1}^d \lambda_{1, j_k}
                         \lesssim       \tau_n^d \exp\!\big( -D_{\exp} \tau_n \| \mathbf{j} \|_1 \big),
                         \text{ where } \|\mathbf{j}\|_1 = \sum_{k=1}^d j_k.
  \end{align}

  We now consider the non-increasing single-index rearrangement $(\lambda_j)_{j \ge 1}$ of the multi-index eigenvalues $(\Lambda_{\mathbf j})_{\mathbf j \in \mathbb{N}^d}$.
	For the counting function
  \begin{align}
    N(q):= \#\{\mathbf{j}\in\mathbb{N}^d : \|\mathbf{j}\|_1 \le q\}, \qquad q \in \mathbb N,
  \end{align}
  it is well known that \( N(q) \asymp q^d \) for \( q \to\infty \).
  Hence, the relation between the rank $j$ and the level $q$ satisfies \( j \asymp q^d \), which is equivalent to \( q \asymp j^{1/d} \).
  Since all the eigenvalues corresponding to multi-indices $\mathbf{j}$ with $\|\mathbf{j}\|_1 \asymp q$ are bounded by $\exp(-D_{\exp}\tau_n q)$, it follows that the $j$-th largest eigenvalue satisfies
  \begin{align}
    \label{eq_Aux_EVDecayForSEKernel}
    \lambda_j \lesssim \tau_n^d \exp\!\big(-D_{\exp} \tau_n q       \big)
              \asymp   \tau_n^d \exp\!\big(-D_{\exp} \tau_n j^{1/d} \big)
              \quad \text{ with } q \asymp j^{1/d},
  \end{align}
  where we potentially adjusted the constant \( D_{\exp} \).

  Next, let \( J = J_n \in \mathbb{N} \) with \( J \le n \) and introduce the notation \( \widehat{\lambda}_j = \widehat{\mu}_{j} / n \) \( j \le n \).
  Then, from Proposition 2 in \cite{WilliamsShaweTaylor2002_StabilityKernelPCA} combined with Markov's inequality, it follows that 
  \begin{align}
           \tr \big( K K_\sigma^{-1} \big)
    & = 
           \sum_{j = 1}^n \frac{\widehat{\mu}_j}{\widehat{\mu}_j + \sigma^2} 
      = 
           \sum_{j = 1}^n \frac{n \widehat{\lambda}_j}{n \widehat{\lambda}_j + \sigma^2} 
      \le 
           J +  \sum_{j = J + 1}^n \frac{n \widehat{\lambda}_j}{n \widehat{\lambda}_j + \sigma^2}.
    \\
    & \le 
           J +  \frac{n}{\sigma^2} \sum_{j = J + 1}^n \widehat{\lambda}_j
      \le 
           J +  \frac{M_n n}{\sigma^2} \sum_{j = J + 1}^n \lambda_j
    \notag
  \end{align}
  with probability converging to one under \( \mathbb{P}_{f_0}^{\otimes n} \) and any \( M_n \to \infty \).
  Restricting to this event and using the inequality in Equation \eqref{eq_Aux_EVDecayForSEKernel}, we obtain
  \begin{align}
             \tr \big( K K_\sigma^{-1} \big) 
    \lesssim 
             M_n \Big( J_n + \frac{n \tau_n^d}{\sigma^2}
                             \sum_{j = J_n + 1}^n \exp \big( -D_{\exp} \tau_n j^{1/d} \big) 
                 \Big).
  \end{align}
  Noting that by Lemma \ref{lem_BoundForTheUpperIncompleteGammaFunction},
  \begin{align}
    & \ \ \ \
              \sum_{j = J_n + 1}^n \exp \big( -D_{\exp} \tau_n j^{1/d} \big)
    \le 
              \int_{J_n}^ \infty \exp \big( -D_{\exp} \tau_n t^{1/d} \big) \, dt
    \\
    & = 
              d \int_{ J_n^{1/d} }^\infty u^{d - 1} \exp \big( -D_{\exp} \tau_n u \big) \, du
    \lesssim 
              \tau_n^{-1} J_n^{1 - 1/d} \exp \Big( -D_{\exp} \tau_n J_n^{1/d} \Big).
    \notag
  \end{align}
  Finally, choosing \( J_n = \tau_n^{-d} \log^d(n) / D_{\exp} \) and adjusting the sequence \( M_n \) to account for constants yields the results.
\end{proof}

\begin{lemma}[Bound for the upper incomplete \( \Gamma \)-function]
  \label{lem_BoundForTheUpperIncompleteGammaFunction}
  For \( p > 0 \), let 
  \begin{align}
    \Gamma(p, x) = \int_x^\infty t^{p - 1} e^{-t} \, dt, \qquad x > 0
  \end{align}
  be the upper incomplete \( \Gamma \)-function.
  Then,
  \begin{align}
    \Gamma(p, x) \le \max \{ 1, 2^{p - 2} \} \big( x^{p - 1} + \Gamma(p) \big) e^{-x}, \qquad x > 0.
  \end{align}
\end{lemma}

\begin{proof}[Proof]
  We have that
  \begin{align}
    \int_x^\infty t^{p - 1} e^{-t} \, dt 
    & = 
    e^{-x} \int_0^\infty (x + s)^{p - 1} e^{-s} \, ds
    \le  
    \max \{ 1, 2^{p - 2} \} e^{-x} \int_0^\infty (x^{p - 1} + s^{p - 1}) e^{-s} \, ds
    \notag
    \\
    & = 
    \max \{ 1, 2^{p - 2} \}  e^{-x} \big( x^{p - 1} + \Gamma(p) \big),
  \end{align}
  where for \( p - 1 \ge 1 \), we have used the convexity estimate \( (x + s)^{p - 1} \le 2^{p - 2} (x^{p - 1} + s^{p - 1}) \) and \( (x + s)^{p - 1} \le x^{p - 1} + s^{p - 1} \) for \( p - 1 < 1 \).
\end{proof}

\begin{lemma}[Number of neighbours bound]
  \label{lem_NumberOfNeighboursBound}
  Let \( X_1, \dots, X_n \sim \) i.i.d. be distributed according to \( G \) on \( \mathcal X \) with Lebesgue density \( g \) which is continuous and strictly positive at \( x \in \mathcal X \).
  For a radius \( r > 0 \), let 
  \begin{align}
    N_i(r) = \sum_{j \ne i} \mathbf{1} \{ \| X_j - X_i \| \le r \}
  \end{align}
  be the number of neighbours of \( X_i \) not further away than \( r \).
  Then, conditionally on \( X_i = x \), with probability at least \( 1 - 2 \exp(- (n - 1) \delta^2 / 4) \), we have that
  \begin{align}
    |N_i(r) - (n - 1) p_r(x)| \le (n - 1) \delta 
  \end{align}
  where \( p_r(x) = \int_{B_r(x)} g(u) d u \asymp g(x) r^{d} \text{Vol}(B_{1}(0)) \) for \( r \to 0 \). 
\end{lemma}

\begin{proof}[Proof]
  For \( r \to 0 \), conditionally on \( X_i = x \in \mathcal X  \), \( N_i(r) \) is distributed as a binomial \( \text{Bin}(n - 1, p_r(x)) \) with 
  \begin{align}
    p_r(x) = \int_{B_r(x)} g(u) d u 
           = \int_{B_1(0)} g(x + r u) r^{d} d u 
           = g(x) r^d\, \text{Vol}(B_1(0)) (1 + o(1)) 
  \end{align}
  for \( r \to 0 \).
  By Chernoff's inequality, c.f. Exercise 2.10 in \cite{BoucheronEtal2013_ConcentrationInequalities}, for \( \delta > 0 \), 
  \begin{align}
    \mathbb P \big( |N_i(r) - (n - 1) p_r(x)| \ge (n - 1) \delta \big) 
    & \le 
    \exp \Big( \frac{- (n - 1) \delta^2}{4} \Big).
  \end{align}
  Note that compared to Exercise 2.10, the bound can be formulated to avoid the appearance of \( p_r(x) \) on the right-hand side.
\end{proof}

\begin{lemma}[Trace term for the Mat\'ern kernel]
  \label{lem_TraceTermForTheMaternKernel}
  Assume that \( \mathcal{X} \subset \mathbb R^d \) is bounded and the design distribution \( G \) on \( \mathcal{X} \)  has a bounded Lebesgue density \( g \).
  Let \( k_{\tau_n} \) be the Mat\'ern kernel with length scale \( \tau_n = n^{-(\nu - \alpha) / (2 \alpha \nu + d \nu)} \) and smoothness \( \nu \).
  Then, for any sequence \( M_n \to \infty \),
  \begin{align}
    \tr K K_ \sigma^{-1} \lesssim M_n n^{d / (2 \alpha + d)}
  \end{align}
  with probability converging to one under \( \mathbb P_{f_0}^{\otimes n} \) for \( n \to \infty \).
\end{lemma}

\begin{proof}[Proof]
  As in Lemma \ref{lem_TraceTermForTheSquaredExponentialKernel}, we have  that for any \( J \in \mathbb N \) and any \( M_n \to \infty \),
  \begin{align}
    \label{eq_Aux_StructuralTraceTermBound}
    \tr K K_\sigma^{-1} \le J + \frac{M_n n}{\sigma^2} \sum_{j = J + 1}^n \lambda_j
  \end{align}
  with probability converging to one under \( \mathbb P_{f_0}^{\otimes n} \).
  
  Now, under the assumptions of Theorem 1 in \cite{Seeger2007_AddendumInformationConsistency} together with the definition of the spectral density of the Mat\'ern prior guarantees that
  \begin{align}
    \lambda_j \lesssim \tau_n^{-2 \nu} j^{-1 - 2 \nu / d}, \qquad j \ge 1.
  \end{align}
  Plugging into the bound in Equation \eqref{eq_Aux_StructuralTraceTermBound} yields
  \begin{align}
    \tr K K_\sigma^{-1} 
    & \lesssim 
    M_n \Big( J + n \tau_n^{- 2 \nu} \sum_{j = J + 1} j^{- 1 - 2 \nu / d}  \Big)
    \lesssim 
    M_n \Big( J + n \tau_n^{- 2 \nu} J^{- 2 \nu / d} \Big).
  \end{align}
  Setting \( J = J_n = n^{d / (2 \alpha + d)} \) finishes the proof.
\end{proof}

\begin{proposition}[Lower bound]
  \label{prp_LowerBound}
  Let \( k \) be a strictly positive definite and bounded kernel, i.e. \( K > 0 \) a.e. and \( \sup k(x, x) < \infty \).
  Further assume that \( f_{0} = \varphi_{1} = \lambda_{1}^{1/2} \phi_{1} \) is the first eigenfunction of the covariance operator \( \Sigma = T_{k} |_{\mathbb{H}} \) of \( k(X_{1}, .) \), see Appendix \ref{ssec_ReproducingKernelHilbertSpaces}.  
  Then, there is an event with probability converging to \( 1 \) on which the approximate posterior is either not well defined in the sense that \(  0 \le K - K \widehat{K}_{\sigma}^{-1} K \le K \) or
  \begin{align}
    \| \mathbb{E}_{\varepsilon} K \widehat{K}_{\sigma}^{-1} Y - f_{0} \|_{n}
    & \ge
    c \sqrt{\lambda_{1}}
    \Big( 
      \Big( \frac{n \lambda_{1}}{2 \| \widehat{K} \|_{\text{op}}} - 1 \Big)
      - 
      \frac{\sqrt{\mathbb{E} \| k(X_{1}, \cdot) \|_{\mathbb{H}}^{4} - \| \Sigma \|_{\text{HS}}^{2}}}
           {(\lambda_{1} - \lambda_{2}) \sqrt{n / \log n}} 
    \Big)^{+}
  \end{align}
  The well definedness condition can be interpreted such that the approximate posterior has a positive semi-definite covariance structure that reduces the variance of the prior.
\end{proposition}

\begin{proof}[Proof]
  \textbf{Step 1: Deterministic analysis of kernel matrix terms.}
  The condition \( K - K \widehat{K}_{\sigma}^{-1} K \le K \) requires that \( K \widehat{K}_{\sigma}^{-1} K \ge 0 \).
  By Sylvester's law of inertia and the fact that \( k \) is strictly positive definite, this is equivalent to \( K^{1/2} \widehat{K}_{\sigma}^{-1} K^{1/2} \ge 0 \), which in turn is equivalent to \( K \widehat{K}_{\sigma}^{-1} \ge 0 \) by the similarity of matrices.
  With the same arguments this is also equivalent to \( \widehat{K}_{\sigma} \ge 0 \).

  With a similar argument \( K - K \widehat{K}_{\sigma}^{-1} K \ge 0 \) is equivalent to \( I - K^{1/2} \widehat{K}_{\sigma}^{-1} K^{1/2} \ge 0 \).
  This implies that that \( K^{1/2} \widehat{K}_{\sigma}^{-1} K^{1/2} \) may only have eigenvalues smaller than 1 and by similarity of the matrices the condition becomes \( K K_{\alpha}^{-1} \le I \).
  In the following we therefore only consider realisations of the design for which \( 0 \le K K_{\sigma}^{-1} \le I \).

  Now, let
  \begin{align}
    S_{n}: \mathbb{H} \to \mathbb{R}^{n}, \qquad h \mapsto (h(X_{1}), \dots, h(X_{n}))^{\top}
  \end{align}
  be the sampling operator and let \( (\widehat{\varphi}_{j})_{j \le n} \) be the eigenfunctions of the empirival covariance operator
  \begin{align}
    \widehat{\Sigma}: 
    =
    \frac{1}{n} \sum_{i = 1}^{n} k(X_{i}, \cdot) \otimes k(X_{i}, \cdot) 
    =
    S_{n}^{*} S_{n}, 
  \end{align}
  where \( S_{n}^{*} \) is the adjoint of \( S_{n} \) with respect to the space \( (\mathbb{R}^{n}, \langle \cdot, \cdot \rangle_{n} ) \).
  Noting that \( S_{n} S_{n}^{*} = n^{-1} K \), the operators \( \widehat{\Sigma} \) and \( n^{-1} K \) have the same eigenvalues \( (\widehat{\lambda}_{j} )_{j \le n} \). 

  We may now write
  \begin{align}
    \| \mathbb{E}_{\varepsilon} K \widehat K_\sigma^{-1} - f_{0} \|_{n}
    & = 
    \| K \widehat{K}_{\sigma}^{-1} S_{n} f_{0} - S_{n} f_{0} \|_{n}
    = 
    \| K \widehat{K}_{\sigma}^{-1} S_{n} \varphi_{1} - S_{n} \varphi_{1} \|_{n}
    \\
    & \ge 
    \| K \widehat{K}_{\sigma}^{-1} S_{n} \widehat{\varphi}_{1} - S_{n} \widehat{\varphi}_{1} \|_{n}
    - 
    \| (K \widehat{K}_{\sigma}^{-1} - I) S_{n} ( \varphi_{1} - \widehat{\varphi}_{1}) \|_{n}.
    \notag
  \end{align}
  We treat the two remaining terms individually.

  Davis-Kahan theorem implies that, by potentially adjusting the sign of \( \widehat{\varphi}_{1} \), 
  \begin{align}
    \| \varphi_{1} - \widehat{\varphi}_{1} \|_{\mathbb{H}} 
    & \le 
    \frac{\| \widehat{\Sigma} - \Sigma \|_{\text{HS}}}{\lambda_{1} - \lambda_{2}}. 
  \end{align}
  Therefore, under the assumption that \( 0 \le K \widehat{K}_{\sigma}^{-1} \le I \), we may estimate
  \begin{align}
    & \ \ \ \
    \| (K \widehat{K}_{\sigma}^{-1} - I) S_{n} ( \varphi_{1} - \widehat{\varphi}_{1}) \|_{n}
    \le 
    \| (K \widehat{K}_{\sigma}^{-1} - I) \|_{\text{op}}
    \|  S_{n} ( \varphi_{1} - \widehat{\varphi}_{1}) \|_{n}
    \\
    & \le 
    C \sqrt{\widehat{\lambda}_{1}} \| \widehat{\varphi}_{1} - \varphi_{1} \|_{\mathbb{H}} 
    \le 
    C \sqrt{\widehat{\lambda}_{1}} 
    \frac{\| \widehat{\Sigma} - \Sigma \|_{\text{HS}}}{\lambda_{1} - \lambda_{2}}.
    \notag
  \end{align}

  For the first term, we use the fact that \( S_{n} \widehat{\varphi}_{1} = \widehat{\lambda}_{1}^{1/2} \widehat{v}_{1} \) where \( \widehat{v}_{1} \) is an eigenvevtor of \( K \) with \( \| \widehat{v}_{1} \|_{n} = 1 \). 
  Therefore, 
  \begin{align}
    & \ \ \ \
    \| K \widehat{K}_{\sigma}^{-1} S_{n} \widehat{\varphi}_{1} - S_{n} \widehat{\varphi}_{1} \|_{n}
    \ge 
    \langle \widehat{v}_{1}, (K \widehat{K}_{\sigma}^{-1} - I) S_{n} \widehat{\varphi}_{1} \rangle_{n}
    = 
    \widehat{\lambda}_{1}^{1/2}
    \langle \widehat{v}_{1}, (K \widehat{K}_{\sigma}^{-1} - I) \widehat{v}_{1} \rangle_{n}
    \\
    & = 
    \widehat{\lambda}_{1}^{1/2} 
    \big( \langle \widehat{v}_{1}, K K_{\sigma}^{-1} \widehat{v}_{1}  \rangle_{n}- 1 \big)
    = 
    \widehat{\lambda}_{1}^{1/2} 
    \big( 
      n \widehat{\lambda}_{1} 
      \langle \widehat{v}_{1}, (I + K_{\sigma}^{-1} E)^{-1} K_{\sigma}^{-1} \widehat{v}_{1} \rangle_{n}
      - 1 
    \big)
    \notag
    \\
    & = 
    \widehat{\lambda}_{1}^{1/2}
    \Big( 
      \frac{n \widehat{\lambda}_{1}}{n \widehat{\lambda}_{1} + \sigma^{2}} 
      \langle \widehat{v}_{1}, (I + K_{\sigma}^{-1} E)^{-1} \widehat{v}_{1} \rangle_{n}
      - 1
    \Big)
    \notag 
    \\
    & = 
    \widehat{\lambda}_{1}^{1/2}
    \Big( 
      \frac{n \widehat{\lambda}_{1}}{n \widehat{\lambda}_{1} + \sigma^{2}} 
      \langle \widehat{v}_{1},
              K_{\sigma}^{-1/2} (I + K_{\sigma}^{-1/2} E K_{\sigma}^{-1/2})^{-1} K_{\sigma}^{1/2}
              \widehat{v}_{1} 
      \rangle_{n}
      - 1
    \Big)
    \notag
    \\
    & = 
    \widehat{\lambda}_{1}^{1/2}
    \Big( 
      \frac{n \widehat{\lambda}_{1}}{n \widehat{\lambda}_{1} + \sigma^{2}} 
      \langle \widehat{v}_{1}, (I + K_{\sigma}^{-1/2} E K_{\sigma}^{-1/2})^{-1} \widehat{v}_{1} \rangle_{n}
      - 1
    \Big).
    \notag
  \end{align}
  We note that \( I + K_{\sigma}^{-1/2} E K_{\sigma}^{-1/2} > 0 \) is equivalent to \( K_{\sigma} + E > 0 \), which is in turn equivalent to \( \widehat{K}_{\sigma} > 0 \).
  We can then use the inequality
  \begin{align}
    \langle v, A v \rangle_n \langle v, A^{-1} v \rangle_n \ge \langle v, v \rangle_n^2,
  \end{align}
  which holds for all \( A > 0 \) and \( v \in \mathbb{R}^{n} \).
  This yields the estimate
  \begin{align}
    & \ \ \ \
    \| K \widehat{K}_{\sigma}^{-1} S_{n} \widehat{\varphi}_{1} - S_{n} \widehat{\varphi}_{1} \|_{n}
    \ge 
    \widehat{\lambda}_{1}^{1/2}
    \Big( 
      \frac{n \widehat{\lambda}_{1}}{n \widehat{\lambda}_{1} + \sigma^{2}} 
      \langle \widehat{v}_{1}, (I + K_{\sigma}^{-1/2} E K_{\sigma}^{-1/2}) \widehat{v}_{1} \rangle_{n}^{-1}
      - 1
    \Big)
    \\
    & = 
    \widehat{\lambda}_{1}^{1/2}
    \Big( 
      \frac{n \widehat{\lambda}_{1}}{n \widehat{\lambda}_{1} + \sigma^{2}} 
      \Big[ 1 + \frac{\langle \widehat{v}_{1}, E \widehat{v}_{1} \rangle_{n}}
                     {n \widehat{\lambda}_{1} + \sigma^{2}} 
      \Big]^{-1} 
      - 1
    \Big)
    \notag
    \\
    & = 
    \widehat{\lambda}_{1}^{1/2}
    \Big( 
      \frac{n \widehat{\lambda}_{1}}{n \widehat{\lambda}_{1} + \sigma^{2}} 
      \Big[ \frac{  n \widehat{\lambda}_{1} + \sigma^{2} 
                  + \langle \widehat{v}_{1}, \widehat{K} \widehat{v}_{1} \rangle_{n}
                 }
                 {n \widehat{\lambda}_{1} + \sigma^{2}} 
      \Big]^{-1} 
      - 1
    \Big)
    \ge 
    \widehat{\lambda}_{1}^{1/2}
    \Big( \frac{n \widehat{\lambda}_{1}}{\| \widehat{K} \|_{\text{op}}} - 1 \Big).
    \notag
  \end{align}

  \textbf{Step 2: Contructing the event.}
  We have
  \begin{align}
    & \ \ \ \
    \mathbb{E} \| \widehat{\Sigma} - \Sigma \|_{\text{HS}}^{2}
    = 
    \mathbb{E} \| \widehat{\Sigma} \|_{\text{HS}}^{2} - \| \Sigma \|_{\text{HS}}^{2}
    = 
    \mathbb{E} \sum_{j = 1}^{n} 
    \langle \widehat{\Sigma} \widehat{\varphi}_{j}, \widehat{\Sigma} \widehat{\varphi}_{j} \rangle
    - 
    \| \Sigma \|_{\text{HS}}^{2}
    \\
    & = 
    \frac{1}{n^{2}}
    \mathbb{E} 
    \Big( 
      \sum_{j = 1}^{n} \sum_{i, i' = 1}^{n} 
      \big\langle \langle k(X_{i},  \cdot), \widehat{\varphi}_{j} \rangle k(X_{i}, \cdot),
                  \langle k(X_{i'}, \cdot), \widehat{\varphi}_{j} \rangle k(X_{i'}, \cdot)
      \big\rangle
    \Big) 
    - 
    \| \Sigma \|_{\text{HS}}^{2}
    \notag
    \\
    & \le 
    \frac{1}{n^{2}} 
    \big( 
      n \mathbb{E} \| k(X_{1}, \cdot) \|_{\mathbb{H}}^{4} 
      + 
      (n^{2} - n) \| \Sigma \|_{\text{HS}}^{2} 
    \big)
    - 
    \| \Sigma \|_{\text{HS}}^{2}
    = 
    \frac{\mathbb{E} \| k(X_{1}, \cdot) \|_{\mathbb{H}}^{4} - \| \Sigma \|_{\text{HS}}^{2}}{n}.
    \notag
  \end{align}
  By an applying Markov inequality, there is an event \( \Omega_{1} \) on which \( \| \widehat{\Sigma} - \Sigma \|_{\text{HS}}^{2} \le (\mathbb{E} \| k(X_{1}, \cdot) \|_{\mathbb{H}}^{4} - \| \Sigma \|_{\text{HS}}^{2}) \log(n) / n \).

  Further, by Proposition 2 in \cite{WilliamsShaweTaylor2002_StabilityKernelPCA} and the fact that \( k \) is bounded, there is an event \( \Omega_{2} \) on which \( \widehat{\lambda}_{1} \ge \lambda_{1} / 2 \) with probability converging to 1.
  The result now follows by restricting to \( \Omega_{1} \cap \Omega_{2} \). 
\end{proof}


\section{Complementary results and examples}
\label{sec_ComplementaryResultsAndExamples}

\subsection{Well definedness of the approximate posterior}
\label{ssec_WellDefinednessOfTheApproximatePosterior}

We show that $ \Psi_{r_n}(\cdot) $ is a well defined Gaussian process whenever \( \| \widehat{K} - K \| < \sigma^2 \).
This follows from the fact that it arises as the posterior of the following perturbed observation model:
\begin{align}
	\label{eq_Complement_PerturbedModel}
	\tilde Y_i = f(X_i) + \xi_i, \qquad i=1,\dots,n,
\end{align}
where the noise vector $\xi = (\xi_i)_{i\le n}$ satisfies \( \xi \sim \mathcal{N}(0,\widehat K - K + \sigma^2 I_n) \) and is independent of the prior process $ W $ on \( f \) (although it may depend on the design $X$).
The condition \(\|\widehat K - K \|_{\rm op} < \sigma^2\) ensures that the covariance matrix $\widehat K - K + \sigma^2 I_n$ is symmetric positive definite, and therefore defines a valid Gaussian noise model.

Let $ \tilde n \in \mathbb N $ and $ \tilde{\mathbf{x}} = (\tilde x_1, \dots, \tilde x_{\tilde n}) \in \mathcal{X}^{\tilde n} $ be some evaluation points.
Under the perturbed model \eqref{eq_Complement_PerturbedModel}, the joint distribution of the observations $\tilde Y$ and the vector $ \tilde{\mathbf{f}} = (f(\tilde x_1), \cdots, f(\tilde x_{\tilde n}))^\top$ is Gaussian,
\begin{align}
        \begin{pmatrix}
          \tilde Y \\ 
          \tilde{\mathbf{f}}
        \end{pmatrix} 
  \sim
        \mathcal{N}\!\left( \begin{pmatrix}
                              0\\
                              0
                            \end{pmatrix},
                            \begin{pmatrix}
                              \widehat K_\sigma     & K(X,        \tilde x) \\
                              K(X, \tilde x)^{\top} & K(\tilde x, \tilde x)
                            \end{pmatrix}
                     \right),
\end{align}
where \( K(\mathbf{x}, \tilde{\mathbf{x}}) = (k(x_i, \tilde x_j))_{i \le n, j \le \tilde n} \).
The assumption $\|\widehat K - K \|_{\mathrm{op}} < \sigma^2$ implies that $\widehat K_\sigma$ is strictly positive definite and therefore invertible.
By the standard formula for conditional distributions of multivariate Gaussian vectors, the conditional distribution of $\mathbf f(\mathbf x)$ given $(X,\tilde Y)$ is
\begin{align}
       \tilde{\mathbf f} \mid X,\tilde Y
	\sim
  \mathcal{N}\!\left( K(X, \tilde{\mathbf x})^{\top} \widehat K_\sigma^{-1} \tilde Y,\;
                      K(\mathbf x,\mathbf x) - K(X,\mathbf x)^{\top}\widehat K_\sigma^{-1}K(X,\mathbf x)
               \right).
\end{align}
Since a Gaussian process is uniquely determined by its finite-dimensional distributions, this shows that $ \Psi_{r_n}(\cdot | (X, Y) = (x, y)) $ is a well-defined Gaussian process and coincides with the posterior distribution \( W | (X, \tilde Y) = (x, y) \) in the perturbed model.
Based on this, we can show that the Kullback-Leibler divergence between the processes \( \Psi_{r_n} \) and \( \Pi_{n}(\cdot | X, Y) \) coincides with the Kullback-Leibler divergence between the processes evaluated at the design points.

\begin{lemma}[Process KL vs. finite dimensional KL]
  \label{lem_ProcessKLVsFiniteDimensionalKL}
  Conditional on the design $ X $ and the observations $ Y $, the KL divergence $ \kl( \Psi_{r_n} , \Pi_{n}( \cdot | X, Y) ) $ between the banded posterior GP and the true posterior GP is equal to the KL divergence $ \kl( \Psi_{r_n}( \mathbf{f} \in \cdot), \Pi_{n}(\mathbf{f} \in \cdot | X, Y)) $ between the finite dimensional Gaussians given by  $ \mathbf{f} = ( f( X_{1} ), \dots f( X_{n} ) )^{\top} $.
\end{lemma}

\begin{proof}[Proof]
  Initially, we introduce some helpful notation.
  Let $ ( \Omega, \mathscr{A}, \mathbb{P} ) $ denote the probability space on which both regression models from Equations \eqref{eq_Intro_GPRegressionModel} and \eqref{eq_Complement_PerturbedModel} are defined conditional on the design $ X $.
  In the following, superscripts will denote push forward measures.
				
  Assuming all quantities exist, we may now write the Kullback-Leibler divergence in terms of Radon-Nikodym derivatives, i.e.,
  \begin{align}
    \label{eq_KLOfGPs}
    & \ \ \ \
        \kl(\Psi_{r_n}(\cdot | Y = y), \Pi_n(\cdot | Y = y)) 
    \notag
    \\
    & = 
        \int \log \frac{d \Psi_{r_n}(\cdot | Y = y)}{d \Pi} \, d \Psi_{r_n}(\cdot | Y = y)
      + 
        \int \log \frac{ d \Pi }{ d \Pi_{n}( \cdot | Y = y) } \, d \Psi_{r_n}(\cdot | Y = y)
    \notag
    \\
    & = 
        \int \log \frac{ d \mathbb{P}^{W | \tilde Y = y} }{ d \mathbb{P}^{W} } 
        \, d \mathbb{P}^{W | \tilde Y = y} 
      + 
        \int \log \frac{ d \mathbb{P}^{W} }{ d \mathbb{P}^{W | Y = y} } 
        \, d \mathbb{P}^{W | \tilde Y = y}
  \end{align}
  and analogously
  \begin{align}
    \label{eq_KLOfFidis}
    & \ \ \ \
        \kl(\Psi_{r_n}(\mathbf f \in \cdot | Y = y), \Pi_n(\mathbf f \in \cdot | Y = y)) 
    \notag
    \\
    & = 
        \int \log \frac{ d \mathbb{P}^{\vec{W} | \tilde Y = y} }{ d \mathbb{P}^{\vec W} } 
             \, d \mathbb{P}^{\vec W | \tilde Y = y} 
      + 
        \int \log \frac{ d \mathbb{P}^{\vec W} }{ d \mathbb{P}^{\vec W | Y = y} } 
             \, d \mathbb{P}^{\vec W | \tilde Y = y},
  \end{align}
  where $ \vec W: = (W(X_1), \dots, W(X_n)) $ denotes the finite dimensional evaluation of the process $ W $. 
  Therefore, it is sufficient to show that the integrands in Equations \eqref{eq_KLOfGPs} and \eqref{eq_KLOfFidis} exist and coincide respectively.
				
  Since any posterior is almost surely absolutely continuous with respect to the prior, we have that the Radon-Nikodym derivatives
  $ d \mathbb{P}^{W |        Y = y} / d \mathbb{P}^{W} $,
  $ d \mathbb{P}^{W | \tilde Y = y} / d \mathbb{P}^{W} $ 
  exist and are well defined.
  Further, via Bayes-theorem,
  \begin{align}
        \frac{ d \mathbb{P}^{W | Y = y} }{ d \mathbb{P}^{W} } 
    & = 
        \frac{d ( p(Y = y | W = f) / p(Y = y) \pi(d f) )}{d \Pi}
      = 
        \frac{p( Y | \vec W = \mathbf f)}{p(Y = y)} 
      >
        0
  \end{align}
  due to the exponential form of the likelihood \( p \).
  Therefore,
  \begin{align}
        \frac{ d \mathbb{P}^{W} }{ d \mathbb{P}^{ W | Y = y} } 
    & = 
        \Big( \frac{ d \mathbb{P}^{W | Y = y} }{ d \mathbb{P}^{W} } \Big)^{-1}
      = 
        \frac{p(Y = y)}{p(Y = y | W)} 
      = 
        \frac{p(Y = y)}{p(Y = y | \vec W)}.
  \end{align}
  We conclude that the right-most terms in Equations \eqref{eq_KLOfGPs} and \eqref{eq_KLOfFidis} exist and coincide.
				
  An analogous argument as above also shows that $ d \mathbb{P}^{W | \tilde Y = y} / d \mathbb{P}^{W} $ only depends on $ \vec W $.
  This finishes the proof, since therefore also the left-most terms in Equation \eqref{eq_KLOfGPs} and \eqref{eq_KLOfFidis} exist and coincide.
\end{proof}


\subsection{Reproducing kernel Hilbert spaces}
\label{ssec_ReproducingKernelHilbertSpaces}

Let \( \mathcal{X} \subset \mathbb R^d \) and let $k: \mathcal{X}^2 \to \mathbb R$ be a continuous positive definite kernel.
The associated reproducing kernel Hilbert space (RKHS) $\mathbb H$ can be characterized in several equivalent ways.\\

\noindent \textbf{Definition and reproducing property.}
The RKHS $\mathbb H$ is defined as the completion of the linear span of functions \( \{k(x,\cdot):x\in\mathcal X\} \), equipped with the canonical inner product
\begin{align}
\left\langle 
\sum_{i=1}^m a_i k(x_i,\cdot),
\sum_{j=1}^n b_j k(s_j,\cdot)
\right\rangle_{\mathbb H}
=
\sum_{i=1}^m\sum_{j=1}^n a_i b_j k(x_i,s_j),
\qquad 
a\in\mathbb R^m,\ b\in\mathbb R^n .
\end{align}
In particular, for every $s\in\mathcal X$, the function $k(s,\cdot)$ belongs to $\mathbb H$ and satisfies the reproducing property
\begin{align}
  f(s) = \langle f,k(s,\cdot)\rangle_{\mathbb H}, \qquad f\in\mathbb H.
\end{align}
Consequently, pointwise evaluation is a continuous linear functional on $\mathbb H$. Moreover, by the Cauchy--Schwarz inequality,
\begin{align}
  |f(x)| \le \|f\|_{\mathbb H}\sqrt{k(x,x)}, \qquad f\in\mathbb H.
\end{align}

\noindent \textbf{Stationary kernels and Bochner's theorem.}
When the kernel is stationary, i.e.
\begin{align}
  k(x, x')=r(x' - x), \qquad x, x' \in \mathcal{X},
\end{align}
for some function $r:\mathbb R^d \to \mathbb R$, a more explicit characterization is available via Fourier analysis.
A continuous function $r: \mathbb R^d \to \mathbb R$ is positive definite if and only if it is the Fourier transform of a finite non-negative measure $\mu$, namely
\begin{align}
  r(h) = \int_{\mathbb R^d} e^{i\langle \omega,h\rangle}\,\mu(d\omega).
\end{align}
If $\mu$ admits a density $S(\omega)$ with respect to Lebesgue measure, then $S$ is called the spectral density of the kernel.
In this case, the associated RKHS consists of the real parts of functions of the form
\begin{align}
  f(t) = \int_{\mathbb R^d} e^{i\langle t,\omega\rangle}\psi(\omega)\,\mu(d\omega), \qquad \psi\in L^2(\mu),
\end{align}
with RKHS norm
\begin{align}
  \|f\|_{\mathbb H}^2 = \inf_{\psi:\,f=\mathcal F\psi} \int_{\mathbb R^d} |\psi(\omega)|^2\,\mu(d\omega).
\end{align}

\noindent \textbf{Covariance operator and spectral decomposition.}
Let
\begin{align}
  (T_kf)(x) = \int_{\mathcal X} k(x,y)f(y)\,d\mathbb G(d y), x \in \mathcal{X}
\end{align}
denote the covariance operator associated with $k$.
Assume that $k$ is measurable, symmetric and positive definite, and that
\begin{align}
  \int_{\mathcal X}k(x,x)\,d\mathbb P_X(x)<\infty.
\end{align}
Then, $T_k$ defines a positive trace-class operator on $L^2(\mathbb P_X)$.
Consequently, there exist eigenpairs $ (\lambda_j,\phi_j)_{j \ge 1} $ satisfying
\begin{align}
  T_k\phi_j=\lambda_j\phi_j, \qquad j \ge 1,
\end{align}
such that
\begin{align}
  k(x,y) = \sum_{j\ge1}\lambda_j\phi_j(x)\phi_j(y), \qquad x, y \in \mathcal{X} 
\end{align}
in $L^2(G \otimes G)$.

Moreover, the associated RKHS admits the characterization
\begin{align}
  \mathbb H = \left\{ f=\sum_{j\ge1}a_j\phi_j: \sum_{j\ge1}\frac{a_j^2}{\lambda_j}<\infty \right\},
\end{align}
with RKHS norm
\begin{align}
  \|f\|_{\mathbb H}^2 = \sum_{j\ge1}\frac{a_j^2}{\lambda_j}.
\end{align}
The decay of the eigenvalues $ (\lambda_j)_{j \ge 1} $ reflects the smoothness of the kernel and the complexity of the associated RKHS. 
Exponential eigenvalue decay, as for squared exponential kernels, corresponds to very smooth functions, whereas polynomial decay, as for Mat\'ern-type kernels, yields Sobolev-type space.\\

\noindent \textbf{Karhunen--Lo\`eve expansion.}
If $W$ is a centered Gaussian process with covariance kernel $k$, then under the above assumptions,
\begin{align}
  W = \sum_{j\ge1} \sqrt{\lambda_j}\,\xi_j\phi_j,
\end{align}
where $(\xi_j)_{j\ge1}$ are i.i.d.\ standard Gaussian random variables, with convergence in $L^2(G)$.
This representation is known as the Karhunen--Lo\`eve expansion.


\subsection{A polynomial kernel}
\label{ssec_APolynomialKernel}

We  consider the stationary polynomial-type kernel of order \(p > 0\), defined as
\begin{align}
  k_{\tau, p}(x, x') = (1 + \tau^{-2} \| x' - x \|^2)^{-p}, \qquad x, x' \in \mathbb{R}^d,
\end{align}
where \( \tau > 0\) is a rescaling parameter.
This kernel induces a Gaussian process with sample paths that are smoother than those of finite-order Mat\'ern processes, but generally less regular than those of the squared exponential kernel. 
Initially, we show that the polynomial kernel is well defined.

\begin{lemma}[Spectral density of the polynomial kernel]
  \label{lem_SpectralDensityOfThePolynomialKernel}
  Let $ p > d/2 $.
  Then, the polynomial kernel \( k_{\tau, p} \)  has the spectral density
  \begin{align}
    f_{\tau, p}(x) = \frac{\Gamma(p - d/2) }{\Gamma(p) (4 \pi)^{d/2}} \tau^d
                     \kappa_{p - d/2, \tau^{-1} \sqrt{2 p - d}}(x),
                     \qquad x \in \mathbb R^d,
  \end{align}
  where 
  \begin{align}
    \kappa_{\nu, \sigma}(x) = \frac{2^{1 - \nu}}{\Gamma(\nu)} 
                              \Big( \frac{\sqrt{2 \nu}}{\sigma} \| x \|_2 \Big)^{\nu} 
                              B_{\nu} \Big( \frac{\sqrt{2 \nu}}{\sigma} \| x \|_2 \Big) ,
                              \qquad x \in \mathbb R^d,
  \end{align}
  is the radial function of the Matérn kernel.
  In particular, the polynomial kernel is a well defined covariance kernel.
\end{lemma}

\begin{proof}
  It is a well known fact that the spectral density of the Matérn process with smoothness \( \nu \) and scaling \( \sigma \) is given by
  \begin{align}
    g_{\nu, \sigma}(x) = \frac{\Gamma(\nu + d/2)}{\pi^{d/2} \Gamma(\nu)} \frac{\sigma^d}{(2 \nu)^{d/2}} 
                         \Big( 1 + \frac{\sigma^2 \| x \|^2}{2 \nu} \Big)^{-(\nu + d/2)},
                         \qquad x \in \mathbb R^d,
  \end{align}
  which, up to constants, is a version of the radial function underlying the polynomial kernel.
  By the Fourier inversion theorem for \( L^2(\mathbb R^d) \cap L^1(\mathbb R^d) \), then
  \begin{align}
    g_{\nu, \sigma}(x) = (2 \pi)^{-d} \int e^{- i \langle x, u \rangle} \kappa_{\nu, \sigma}(u) d u
                       = (2 \pi)^{-d} \int e^{  i \langle x, u \rangle} \kappa_{\nu, \sigma}(u) d u, 
                       \qquad x \in \mathbb R^d.
  \end{align}
  By Bochner's theorem, it then follows that the polynomial kernel is a well defined covariance kernel.
  The explicit formula for the spectral density then follows from setting \( \sigma = \sqrt{2 \nu} / \tau \) and \( \nu = p - d/2 \).
\end{proof}

Similarly to the squared exponential case, using a fixed scale parameter \( \tau \) leads to suboptimal performance when estimating a regression function \(f_0 \in C^\alpha(\mathbb{R}^d)\) due to a mismatch between the prior regularity and that of \( f_0 \).
To address this issue, we set \( \tau = \tau_n \), which shrinks the kernel as \( n \) increases.
In the next Lemma, we show that this kernel fits in the framework of \cite{VanDerVaartVanZanten2009_AdaptiveBayesianEstimation} which yields the optimal choice \( \tau_n = n^{-1/(2 \alpha + d)}\). 

\begin{lemma}[Concentration function inequality for the polynomial kernel]
  \label{lem_ConcentrationFunctionInequalityForThePolynomialKernel}
  Let \( \mathcal X = [0, 1]^d \), \( f_0 \in C^\alpha([0, 1]^d) \) and \( p = d/2 + \ell + 1/2 \) for some \( \ell \in \mathbb N \).
  Then, the polynomial kernel \( k_{\tau_n,  p} \)  with lengthscale \( \tau_n = n^{- 1 / (2 \alpha + d)} \) satisfies the concentration function inequality in Equation \eqref{eq_Main_ConcentrationFunction} with \( \eta_n = n^{- \alpha / (2 \alpha + d)} \log^{1 + d/2} n \).
\end{lemma}

\begin{proof}[Proof]
  By Lemma \ref{lem_SpectralDensityOfThePolynomialKernel}, the spectral density of the 
  \begin{align}
    & \ \ \ \
    f_{\tau, p}(x)   = \frac{\Gamma(p - d/2) }{\Gamma(p) (4 \pi)^{d/2}} \tau_n^d
                       \kappa_{p - d/2, \tau_n^{-1} \sqrt{2 p - d}}(x),
                       \\
                   & = \frac{\Gamma(\ell + 1/2)}{\Gamma(d/2 + \ell + 1/2) (4 \pi)^{d/2}} \tau_n^d
                       \kappa_{\ell + 1/2, \tau_n^{-1} \sqrt{2 \ell + 1}}(x)
                       \notag
                       \\
                   & = \frac{\Gamma(\ell + 1/2)}{\Gamma(d/2 + \ell + 1/2) (4 \pi)^{d/2}} 
                       \frac{\Gamma(\ell + 1)}{\Gamma(2 \ell + 1)} 
                       \tau_n^d
                       \notag
                       \exp(-\tau_n \| x \|_2) 
                       \sum_{q = 0}^{\ell} 
                       \frac{(\ell + q)!}{q! (\ell - q)!}  
                       \Big( \frac{\sqrt{8 \ell + 4} \tau_n \| x \|_2}{\sqrt{2 \ell + 1}} \Big)^{\ell - q},
                       \notag
  \end{align}
  where we have used the fact that for \( p = d/2 + \ell + 1/2 \), the radial function of the Mat\'ern kernel has an explicit form.

  It is easy to check that setting \( \tau = \tau_n = 1 \) and \( \delta < 1 \),
  \begin{align}
    \int e^{\delta \| x \|_2} f_{1, p}(x) \, d x < \infty,
  \end{align}
  which is condition (3.3) in \cite{VanDerVaartVanZanten2009_AdaptiveBayesianEstimation}.
  Consequently, Lemmas 4.3 and 4.6 from the same paper yield the inequalities
  \begin{align}
    \tau_n^{-d} \log^{1 + d}(1/ (\tau_n \eta_n)) \lesssim n \eta_n^2,
    \qquad 
    \tau_n^{- d} \lesssim n \eta_n^2
    \qquad \text{ and } \qquad 
    \tau_n^{\alpha} \lesssim \eta_n,
  \end{align}
  for the concentration function inequality \eqref{eq_Main_ConcentrationFunctionInequality} to be satisfied with rate \( \eta_n \).
  This is true, however, for the choices \( \tau_n = n^{-1 / (2 \alpha + d)} \) and the rate \( \eta_n = n^{- \alpha / (2 \alpha + d)} \log^{1 +d/2} n \), which finishes the proof.
\end{proof}

In similar fashion to the examples in Section \ref{sec_Examples}, we obtain a contraction rate result for the approximate posterior.
Note, however, that in this case, the required radius for the contraction regime does not coincide with the one required for the Kullback-Leibler divergence to converge to zero.

\begin{proposition}[Result for the polynomial kernel]
  \label{prp_ResultForThePolynomialKernel}
  In model \eqref{eq_Intro_GPRegressionModel}, assume that \( \mathcal X = [0, 1]^d \) and the desgin observations are distributed according to \( G \) on \( \mathcal{X} \) with bounded Lebesgue density \( g \).
  Let \( f_0 \in C^{\alpha}(\mathcal X) \) and let \( k_{\tau_n, p} \) be the rescaled polynomial kernel with lengthscale \( \tau_n = n^{-1 / (2 \alpha + d)} \) and \( p > d/2 + 2 \alpha \).
  Then, the approximate posterior with truncation radius \( r_n = n^{-(p - 2 \alpha - d/2) / ((2 \alpha + d) p)} \) 
  satisfies that for arbitrary \( M_n \to \infty \)
  \begin{align}
    \Psi_{r_n}(d_H(\cdot, f_0) \ge M_n n^{- \alpha / (2 \alpha + d)}) \to 0
  \end{align}
  in probability with respect to \( \mathbb P_{f_0}^{\otimes n} \) for \( n \to \infty \).
  Further, the number of non-zero entries per row of \( \widehat K \) is at most of order \( m_n = n r_n^d = n^{((2 \alpha + d/2) d/p + 2 \alpha) / (2 \alpha + d)} \) with probability converging to one.

  Moreover, if \( p > d + 2 \alpha \) and \( r_n = n^{-(p - 2 \alpha - d) / ((2 \alpha + d) p)} \log^{1/(2p)} n \), then also
  \begin{align}
    \kl(\Psi_{r_n}, \Pi_n(\cdot, | X, Y)) \to 0
  \end{align}
  in probability with respect to \( \mathbb P_{f_0}^{\otimes n} \) for \( n \to \infty \).
  The number of non-zero entries per row of \( \widehat K \) then is at most of order \( m_n = n r_n^d = n^{((2 \alpha + d) d/p + 2 \alpha) / (2 \alpha + d)} \) with probability converging to one.
\end{proposition}

\begin{proof}[Proof]
  We initially consider the first statement.
  As in the examples before, we have the estimate
  \begin{align}
    \| \widehat K - K \|_{\text{op}} 
    & \le 
    \max_{i \le n} \sum_{X_j: \| X_j - X_i \|_2 > r_n} k_{\tau_n, p}(X_i, X_j) 
    \\
    & \le 
    \max_{i \le n} \sum_{X_j: \| X_j - X_i \|_2 > r_n} \Big( 1 + \frac{r_n^2}{\tau_n^2} \Big)^{- p}
    \lesssim 
    n \Big( \frac{r_n}{\tau_n} \Big)^{-2p}
    \notag
    \\
    & = 
    n^{-2 \alpha / (2 \alpha + d)}
    \to 
    0.
    \notag
  \end{align}
  In particular, \( \| \widehat K - K \|_{\text{op}} < \sigma^2 / 2 \) eventually for all \( n \in \mathbb N \).
  Further, by Theorem \ref{thm_KullbackLeiblerBound}, 
  \begin{align}
    \kl(\Psi_{r_n}, \Pi_n(\cdot | X, Y)) 
    & \lesssim 
    M_n \| \widehat K - K \|_{\text{op}} (\tr K K_\sigma^{-1} + n \eta_n^2) 
    \\
    & \lesssim
    M_n n n^{-2 \alpha / (2 \alpha + d)}
    = 
    M_n n \eta_n^2.
    \notag
  \end{align}
  As in the proof of Theorem \ref{thm_ContractionRates}, we now obtain the contraction rate \( \eta_n \) around \( f_0 \).
  For the statement about the non-zero row entries, we consider again Lemma \ref{lem_NumberOfNeighboursBound} with \( n r_n^d \asymp n^{((2 \alpha + d/2) d/p + 2 \alpha) / (2 \alpha + d)} \) and \( \delta = \| g \|_\infty r_n^d \text{Vol}(B_1(0)) \). 

  The second statement follows from the same bounds with the new choice of \( r_n \).
\end{proof}

\begin{remark}[An improved choice of the radius]
  \label{rem_AnImprovedChoiceOfTheRadius}
  In the proof of Proposition \ref{prp_ResultForThePolynomialKernel}, we have used the straigt forward estimate \( \tr K K_\sigma^{-1} \le n \).
  As in the other examples, the trace term should also satisfy the bound \( \tr K K_\sigma^{-1} \lesssim n \eta_n^2 \).
  However, proving the required bounds for the population eigenvalues goes beyond the scope of this paper.
  Note that the result from \cite{Seeger2007_AddendumInformationConsistency} that was used in Lemma \ref{lem_TraceTermForTheMaternKernel} does not apply here.

  For completeness, we note that under the assumption \( \tr K K_\sigma^{-1} \lesssim n \eta_n^2 \), the required radius for contraction could be improved to 
  \( r_n = n^{-(p - \alpha - d/2) / ((2 \alpha + d) p)} \log^{1/(2p)} n  \) resulting in \( m_n = n r_n^d = n^{((\alpha + d/2) d/p + 2 \alpha) / (2 \alpha + d)} \log^{d / (2p)} n \) of non-zero row entries.
  Analogously, for the convergence of the Kullback-Leibler divergence, we would need \( r_n = n^{-(p - \alpha - d) / ((2 \alpha + d) p)} \\ \log^{1/(2p)} n \) resulting in \( m_n = n r_n^d = n^{((\alpha + d) d/p + 2 \alpha) / (2 \alpha + d)} \log^{d / (2p)} n \) of non-zero row entries.
\end{remark}


\subsection{Additional numerical results}
\label{ssec_AdditionalNumericalResults}

We consider one-dimensional nonparametric regression with observations
\begin{align}
  Y_i=f_0(X_i)+\sigma\varepsilon_i,\qquad i=1,\ldots,n,
\end{align}
where the design points $X_i$ and errors $\varepsilon_i$ are mutually independent standard normal random variables. 
The noise standard deviation is $\sigma=0.2$.
The true regression function is  given by
\begin{align}
  f_0(x)=|x-0.4|^{\beta}-|x-0.2|^{\beta}, \qquad \beta=0.6.
\end{align}
This function has a smoothness equal to $\beta=0.6$. 
We use the rescaled squared exponential kernel as a prior.

We compare the full Gaussian process  posterior mean with the banded approximation, the conjugate gradient (CG) approximation, and distributed Method for $n\in\{10\,000,20\,000, \\ 30\,000\}$.
Each configuration is evaluated over $R=100$ repetitions.
For a test grid $x_1,\ldots,x_L$, the error in repetition $b$ is
\begin{align}
  \mathrm{MSE}^{(b)} = \frac{1}{L}\sum_{\ell=1}^{L} \bigl(\widehat f^{(b)}(x_\ell)-f_0(x_\ell)\bigr)^2.
\end{align}
We report its average and empirical standard deviation across repetitions.

The banding radii are $r\in\{0.75r_n,r_n,1.25r_n\}$, where $r_n=\tau_n\sqrt{2\log n}$ and $\tau_n=n^{-1/2.2}$.
For distributed Method, we consider the number of local machines $m\in\{2,5,50\}$.
For CG, we compare a maximum budgets of $100$ and $150$ iterations with the reference budgets $q_n=303,447,559$, respectively.

Fitting and prediction times are measured separately, and their sum defines total runtime.  the average paired speedups
\begin{align}
  S_{\mathrm{fit}} = \frac{1}{R}\sum_{b=1}^{R} \frac{T_{\mathrm{full,fit}}^{(b)}} {T_{\mathrm{method,fit}}^{(b)}}, \qquad S_{\mathrm{total}} 
                   = \frac{1}{R}\sum_{b=1}^{R} \frac{T_{\mathrm{full,total}}^{(b)}} {T_{\mathrm{method,total}}^{(b)}}.
\end{align}
$T_{\mathrm{fit}}$ is the time needed to calculate the posterior mean and $T_{\mathrm{prediction}}$is the time needed to calculate the mean in new $2000$ test points.
A speedup larger than one therefore indicates a computational gain relative to the full GP.

\begin{table}[H]
\centering
\small
\setlength{\tabcolsep}{4pt}
\begin{tabular}{lrrrr}
\hline
Method and setting
& Mean MSE (SD)
& $S_{\mathrm{fit}}$
& $S_{\mathrm{total}}$ \\
\hline

\multicolumn{4}{l}{$n=10\,000$} \\
Full GP
& $2.610\;(0.315)$  & 1.00 & 1.00 \\
Banded ($r=0.75r_n$)
& $2.610\;(0.315)$  & 5.53 & 4.92 \\
Banded ($r=r_n$)
& $2.610\;(0.315)$  & 3.37 & 3.17 \\
Banded ($r=1.25r_n$)
& $2.610\;(0.315)$  & 2.25 & 2.18 \\
CG ($q_{\max}=100$)
& $2.642 \;(0.315)$  & 1.10 & 1.11 \\
CG ($q_{\max}=150$)
& $2.611\;(0.315)$  & 0.79 & 0.80 \\
CG ($q_{\max}=303$)
& $2.610\;(0.315)$  & 0.41 & 0.42 \\
Distributed  ($m=2$)
& $2.670\;(0.293)$  & 2.38 & 2.29 \\
Distributed  ($m=5$)
& $3.091\;(0.353)$  & 10.48 & 8.24 \\
Distributed  ($m=50$)
& $6.301\;(0.706)$  & 142.13 & 28.61 \\
\hline

\multicolumn{4}{l}{$n=20\,000$} \\
Full GP
& $1.858\;(0.209)$  & 1.00 & 1.00 \\
Banded ($r=0.75r_n$)
& $1.858\;(0.209)$  & 9.16 & 8.54 \\
Banded ($r=r_n$)
& $1.858\;(0.209)$  & 5.49 & 5.28 \\
Banded ($r=1.25r_n$)
& $1.858\;(0.209)$  & 3.66 & 3.58 \\
CG ($q=100$)
& $1.968\; (0.215)$  & 1.87 & 11.87 \\
CG ($q=150$)
& $1.863\;(0.209)$  & 1.33 & 1.33 \\
CG ($q=447$)
& $1.858\;(0.209)$  & 0.53 & 0.53 \\
Distributed  ($m=2$)
& $1.899\;(0.191)$  & 3.24 & 3.18 \\
Distributed  ($m=5$)
& $2.144\;(0.202)$  & 12.63 & 11.44 \\
Distributed  ($m=50$)
& $6.174\;(0.557)$  & 307.06 & 81.87 \\
\hline

\multicolumn{4}{l}{$n=30\,000$} \\
Full GP
& $1.516\;(0.168)$  & 1.00 & 1.00 \\
Banded ($r=0.75r_n$)
& $1.516\;(0.168)$  & 13.21 & 12.57 \\
Banded ($r=r_n$)
& $1.516\;(0.168)$  & 7.81 & 7.59 \\
Banded ($r=1.25r_n$)
& $1.516\;(0.168)$  & 5.07 & 4.98 \\
CG ($q=100$)
& $1.727\; (0.194)$  & 2.63 & 2.62 \\
CG ($q=150$)
& $1.525\;(0.166)$  & 1.88 & 1.87 \\
CG ($q=559$)
& $1.516\;(0.168)$  & 0.67 & 0.67 \\
Distributed  ($m=2$)
& $1.569\;(0.154)$  & 2.94 & 2.92 \\
Distributed  ($m=5$)
& $1.763\;(0.152)$  & 15.07 & 14.22 \\
Distributed  ($m=50$)
& $5.556\;(0.392)$  & 403.85 & 149.39 \\
\hline
\end{tabular}
\caption{
Accuracy and computational performance over $100$ repetitions.
Both mean MSE and its standard deviation, reported in parentheses,
are multiplied by $10^3$.
The CG parameter $q$ denotes the  iteration budget;
iterations may stop earlier.
}
\end{table}

 
   
\end{appendices}

\ \\

\textbf{Funding and acknowledgements.}
We thank the anonymous reviewers for their helpful comments.
Co-funded by the European Union (ERC, BigBayesUQ, project number: 101041064).
Views and opinions expressed are however those of the author(s) only and do not necessarily reflect those of the European Union or the European Research Council.
Neither the European Union nor the granting authority can be held responsible for them.
The research of Bernhard Stankewitz has been partially funded by the Deutsche Forschungsgemeinschaft
(DFG) – Project-ID 318763901-SFB1294.

\printbibliography

\end{document}